\documentclass[12pt,reqno]{article}
\usepackage{amsmath,amsthm,amsfonts,amssymb,amscd,amsbsy}
\usepackage{graphics,color}
\usepackage{url}

\newcommand{\ov}{\overline}

\begin{document}

\newcommand{\nb}{{\mathbb N}}
\newcommand{\la}{{\lambda}}
\newcommand{\al}{{\alpha}}
\newcommand{\be}{{\beta}}
\newcommand{\zb}{{\mathbb Z}}
\newcommand{\qb}{{\mathbb Q}}
\newcommand{\re}{{\mathbb R}}
\newcommand{\te}{{\theta}}

\newcommand{\ga}{{\gamma}}
\newcommand{\ve}{{\varepsilon}}
\newcommand{\vr}{{\varphi}}
\newcommand{\und}{\underline}
\def\ovt{\underline{\theta}}

\centerline{\huge\bf Geometric properties of the Markov}

\medskip

\centerline{\huge\bf and Lagrange spectra}

\vskip .3in

\centerline{\sc Carlos Gustavo Moreira}

\vskip .2in

\centerline{IMPA}
\vskip .2 in
{\small \hskip 2.2in \it Dedicated to Jean-Christophe Yoccoz and Welington de Melo}
\vskip .2in

{\bf Abstract:} We prove several results on (fractal) geometric properties of the classical Markov and Lagrange spectra. In particular, we prove that the Hausdorff dimensions of intersections of both spectra with half-lines always coincide, and may assume any real value in the interval $[0, 1]$.

\section {Introduction and statement of the results}

\vskip .1in

Let $\al$ be an irrational number. According to Dirichlet's theorem, the
inequality \linebreak $|\al-\frac pq|<\frac1{q^2}$ has infinitely many
rational solutions $\frac pq$. Hurwitz improved this result by proving that
$|\al-\frac pq|<\frac1{\sqrt 5 q^2}$ also has infinitely many rational
solutions $\frac pq$ for any irrational $\al$, and that $\sqrt 5$ is the largest
constant that works for any irrational $\al$. However, for particular values of
$\al$ we can improve this constant.

More precisely, if we define $k(\al):=\sup\{k>0\mid  |\al-\frac pq|<\frac 1{
kq^2}$ has infinitely many rational solutions $\frac pq\}=
\limsup_{q\to+\infty}\, ((q|q \al-p|)^{-1})$, we have
$k(\al)\ge\sqrt 5$, $\forall \al\in\re\setminus \Bbb Q$ and $k\left(
\frac{1+\sqrt 5}2 \right)=\sqrt 5$. 
\vskip .1in
\noindent {\bf Definition 1:} The {\it Lagrange spectrum} is the set $L=\{k(\al) \mid
\al\in\re\setminus \Bbb Q$, $k(\al)<+\infty\}$.
\vskip.1in
Hurwitz-Markov theorem determines the first
element of $L$, which is $\sqrt 5$. This set $L$ encodes many diophantine properties of real
numbers. The study of the geometric structure of $L$ is a classical subject, which began with
Markov, proving in 1879 ([Ma]) that
$$
L\cap(-\infty, 3)=\large\{k_1=\sqrt 5< k_2=2 \sqrt 2<k_3= \frac{\sqrt{221}}5<\cdots\large\}
$$
where $k_n$ is a sequence (of irrational numbers whose squares are rational) converging to $3$ - more precisely, the elements $k_n$ of $L\cap(-\infty, 3)$ are the numbers the form $\sqrt{9-\frac4{z^2}}$, where $z$ is a positive integer such that there are other positive integers $x, y$ with $x\le y\le z$ and $x^2+y^2+z^2=3xyz$. This means that the ``beginning'' of the set $L$ is discrete. 

This is not true for the whole set $L$. Indeed, M.
Hall proved in 1947 ([H]) a result on sums of continued fractions with coefficients bounded by $4$: if $C_4=\{\al=[0; a_1, a_2,\ldots]\in [0,1]| a_n\le 4, \forall n\ge 1\}$ then $C_4+C_4=\{x+y|x, y\in C_4\}=[\sqrt{2}-1, 4(\sqrt{2}-1)]$. This implies that $L$ contains the whole half line $[6,+\infty)$.
In 1975, G. Freiman ([F]) determined the biggest half line $[c, +\infty)$ that is contained in
$L$: he proved that
$$
c=\frac{2221564096 + 283748\sqrt{462}}{491993569} \cong 4,52782956616\dots\,.
$$
This half-line is known as {\it Hall's ray}.

These last two results are based on the study of arithmetic sums of regular Cantor
sets, whose relationship with the Lagrange spectrum will be explained below.

Since the best rational approximations of an irrational number are its convergents (from its continued fraction representation), it is not surprising that $k(\al)$ is related to the continued fraction of $\al$. In fact, if the continued fraction of $\al$ is 
$$
\al = [a_0;a_1,a_2,\dots] = a_0 +\cfrac1{a_1 + \cfrac1{a_2 + {\atop\ddots}}}
$$
then we have the
following formula:
$$
k(\al)=\limsup_{n\to\infty}(\al_n+\be_n),
$$
where $\; \al_n=[a_n;a_{n+1},a_{n+2},\dots]\;\text{ and}\;\be_n=[0; a_{n-1}, a_{n-2},\dots,a_1]$.
This follows from the equality
$$
 \left|\al-\frac{p_n}{q_n}\right| =\frac1{q_n(\al_{n+1}q_n+q_{n-1})}=\frac1{(\al_{n+1}+\be_{n+1}) q_n^2}\: ,\quad \forall
n\in\Bbb N.
$$

This formula for $k(\al)$ implies the following alternative
definition of the Lagrange spectrum $L$, due to Perron ([P]): let $\Sigma=({\Bbb N^*})^{\zb}$ be the set of all bi-infinite sequences of
positive integers. If $\und\theta=(a_n)_{n\in\zb}\in \Sigma$, let $\al_n=[a_n;a_{n+1},a_{n+2},\dots]$ and $\be_n=[0;a_{n-1},a_{n-2},\dots], \forall n \in \zb$. We define
$$f(\und\theta)=\al_0+\be_0=[a_0; a_1,a_2,\dots]+[0; a_{-1}, a_{-2},\dots].$$
Then, if $\sigma\colon\Sigma\to\Sigma$ is the shift map defined by
$\sigma((a_n)_{n\in\zb})=(a_{n+1})_{n\in\zb}$, then the Lagrange spectrum is equal to $L=\{\limsup_{n\to\infty} f(\sigma^n \und\theta), \und\theta\in\Sigma\}$.

In this context we can also define the {\it Markov spectrum}.
\vskip .1in
\noindent {\bf Definition 2:} The Markov spectrum is the set 
$M=\{\sup_{n\in\zb} f(\sigma^n\und\theta), \und\theta\in\Sigma\}$. 
\vskip .1in
It also has an arithmetical interpretation (see [P]), namely
$$
M=\{(\inf_{(x,y)\in{\zb}^2\setminus(0,0)}  |f(x,y)|)^{-1},\quad
f(x,y)=a x^2 + bxy+cy^2, \quad b^2-4ac=1\}.
$$
It is well-known (see [CF]) that $M$ and $L$ are closed sets of the real line
and $L\subset M$. In particular, $M$ also contains the Hall's ray $[c,+\infty)$. Freiman also proved in [F] that this is the biggest half-line contained in $M$. 

In this paper, we study the geometrical behaviour of $L$ and $M$ between $3$ and $c$. Consider the function $d:\re\to [0,1]$ defined by $d(t)=HD(L\cap(-\infty,t))$, where $HD$ denotes Hausdorff dimension (see [Fa] for the definitions and basic properties of the notions
of dimension used in this paper). We will prove the following results about the Markov and Lagrange spectra:

\bigskip

{\bf Theorem 1:} {\sl Given $t\in\re$ we have $d(t)=HD(L\cap(-\infty,t))=HD(M\cap(-\infty,t))=\ov{\dim}(L\cap(-\infty,t))=\ov{\dim}(M\cap(-\infty,t))$, where $\ov{\dim}$ denotes upper box dimension. Moreover, $d(t)$ is a continuous non-decreasing surjective function from $\re$ to $[0,1]$, and we have:

i) $d(t)=\min\{1,2 D(t)\}$, where
$D(t):=HD(k^{-1}(-\infty,t)) = HD(k^{-1}(-\infty,t])$ is a continuous function from $\re$ to $[0,1)$.

ii) $\max\{t\in\re\mid d(t)=0\}=3$.

iii) There is $\delta>0$ such that $d(\sqrt{12}-\delta)=1$.
}

\bigskip

This theorem solves affirmatively Problem 3 of [B]. It also gives some answers to Problem 5 of the same paper: the continuous function $d(t)=HD(L\cap(-\infty,t))$, which coincides (for $t>0$) with $\sigma(1/t)$, in the notation of [B], is a Cantor stair function: it is constant in the connected components of the complement of $L \cap (-\infty, t_1]$, where $t_1:=\min \{t \in \re \mid d(t)=1\}\le \sqrt{12}-\delta<\sqrt{12}$; notice that $L \cap (-\infty, t_1]$ is a compact set with zero Lebesgue measure, and so with empty interior. On the other hand, we have the following

{\bf Corollary:} $d(t)$ is not a  H\"older continuous function.

{\bf Proof:}  Suppose by contradiction that $d(t)$ is H\"older continuous with exponent $\alpha>0$. By the previous theorem, there is $\ve>0$ such that $0<d(3+\ve)<\alpha$. Since $d(t)$ is constant in connected components of the open set $\mathbb R \setminus L$, the function $d$ maps the set $L \cap (-\infty,3+\ve]$, whose Hausdorff dimension is $d(3+\ve)<\alpha$, to the nontrivial interval $[0,d(3+\ve)]$. This is a contradiction, since the image of any set of Hausdorff dimension smaller than $\alpha$ by a H\"older continuous function with exponent $\al$ has zero Lebesgue measure (and indeed Hausdorff dimension smaller than one).
\qed

{\bf Remark:} The proof of Theorem 1 doesn't give any estimate on the modulus of continuity of $d(t)$. However it is possible to give such an estimate by modifying the proof. See the discussion at the end of section 6. 

The proof of Theorem 1 is based on the idea of approximating parts of the
spectra from inside and from outside by sums of regular Cantor sets. Theorem 1 uses
techniques developed in a joint work with J.C. Yoccoz about sums of Cantor
sets that implies that the sum of two non essentially affine regular Cantor
sets have Hausdorff dimension equal to the minimum between one and the sum of
their Hausdorff dimensions. This result will be discussed in the next section.
The other results are Theorems 2 and 3 below.

Bugeaud defines in [B], for $c>0$, the sets
$$\text{Exact}(c)=\{\alpha\in \re| \left|\alpha-\frac{p}{q}\right|<\frac{c}{q^2}\text{ for infinitely many }(p, q)\in \zb \times \nb^*\text { but, for every } \ve>0,$$
$$\left|\alpha-\frac{p}{q}\right|<\frac{c-\ve}{q^2}\text{ has only a finite number of solutions }(p, q)\in \zb \times \nb^*\} \text{ and}$$
$$\text{Exact}'(c)=\{\alpha\in \re|\text{ For every } \ve>0, \left|\alpha-\frac{p}{q}\right|<\frac{c+\ve}{q^2}\text{ for infinitely many }(p, q)\in \zb \times \nb^*\text { but}$$
$$\left|\alpha-\frac{p}{q}\right|<\frac{c}{q^2}\text{ has only a finite number of solutions }(p, q)\in \zb \times \nb^*\}.$$
Clearly Exact$(c) \,\cup\,$Exact$'(c)=k^{-1}(c^{-1})$.

{\bf Theorem 2:} {\sl $\lim_{c\to 0}HD(\text{Exact}(c))=\lim_{c\to 0}HD(\text{Exact}'(c))=1$. Consequently, $\lim_{t\to\infty}HD(k^{-1}(t))=1$ and $\lim_{t\to\infty}HD(k^{-1}(-\infty,t))=\lim_{t\to\infty}D(t)=1$.}

This solves affirmatively Problem 4 of [B].

We also prove a result on the topological structure of the Lagrange spectrum $L$.

{\bf Theorem 3:} {\sl $L'$ is a perfect set, i.e., $L''=L'$.}

The proof of this theorem uses the fact that an element of the Lagrange spectrum associated to an infinite sequence $\und\theta$ is accumulated by infinitely many sums of the type $\al_n+\be_n$, which is not necessarily true for elements of the Markov spectrum. The question of whether $M''=M'$ is still open.

There are still some important questions left on the structure of the Markov and Lagrange spectra. For instance: 

{\bf 1)} Consider the function $d_{loc}: L'\to \re$ given by $d_{loc}(t)=\lim_{\ve \to 0} HD(L\cap (t-\ve, t+\ve))$. Is $d_{loc}$ a non-decreasing function ?

{\bf 2)} Describe the geometric structure of the difference set $M\setminus L$.

{\bf 3)} Let, as before, $t_1=\min \{t \in \re \mid d(t)=1\}<\sqrt{12}$. Is it true that \break$t_1=\inf$ int $L$?
This would imply that int$(L \cap (-\infty, \sqrt{12}])\neq\emptyset$, which in its turn implies that int$(C_2+C_2)\neq\emptyset$, where $C_2=\{\al=[0; a_1, a_2,\ldots]\in [0,1]| a_n\le 2, \forall n\ge 1\}$.

We should mention that, in relation to question 2), there are some progresses in recent preprints by the author and C. Matheus: in [MM3] and [MM1], the authors describe the intersections of $M\setminus L$ with the maximal intervals not intersecting $L$ containing the first examples $\gamma=3.11812017815993...$ and $\alpha_{\infty}=3.293044265...$ by Freiman of elements in $M\setminus L$. These intersections have positive Hausdorff dimensions (which coincide with the Hausdorff dimensions of certain regular Cantor sets we describe in these works). In [MM4], the authors exhibit a regular Cantor set diffeomorphic to $C_2$ (and so with Hausdorff dimensions larger than $0.53128$) near $3.7096998597502$ contained in $M\setminus L$. And, in [MM3], they prove that the Hausdorff dimension of $M\setminus L$ is smaller than $0.986927$, and indicate, using heuristic estimates, that this Hausdorff dimension is smaller than $0.888$. 

In relation to question 4), the question whether int$(C_2+C_2)\neq\emptyset$ was posed in page 71 of [CF].

I would like to thank Yann Bugeaud, Aline Gomes Cerqueira, Carlos Matheus, T\'ulio Carvalho and Yuri Lima for helpful comments and suggestions which substantially improved this work.

\section{A dimension formula for arithmetic sums of regular Cantor sets}

We say that $K \subset \re$ is a {\it regular Cantor set of class $C^k$}, $k\ge 1$, if:

\begin{itemize}

\item[i)] there are disjoint compact intervals $I_1,I_2,\dots,I_r$ such that $K\subset I_1 \cup \cdots\cup I_r$ and the boundary of each $I_j$ is contained in $K$;

\item[ii)] there is a $C^k$ expanding map $\psi$ defined in a neighbourhood of $I_1\cup I_2\cup\cdots\cup I_r$ such that $\psi(I_j)$ is the convex hull of a finite union of some intervals $I_s$ satisfying:

\begin{itemize}

\item[ii.1)] for each $j$, $1\le j\le r$, and $n$ sufficiently big, $\psi^n(K\cap I_j)=K$;

\item[ii.2)] $K=\bigcap\limits_{n\in\mathbb N} \psi^{-n}(I_1\cup I_2\cup\cdots\cup I_r)$.

\end{itemize}
\end{itemize}

We say that $\{I_1,I_2,\dots,I_r\}$ is a {\it Markov partition} for $K$ and that $K$ is {\it defined} by $\psi$.

Let $K$ be regular Cantor sets of class $C^2$ defined by the expansive function $\psi$. It is a general fact, due originally to Poincar\'e, that, given a periodic point $p$ of period $r$ of $\psi$, there is a $C^2$ diffeomorphism $h$ of the support interval $I$ of $K$ such that $\tilde\psi=h^{-1}\circ \psi^r\circ h$ is affine in $h^{-1}(J)$, where $J$ is the connected component of the domain of $\psi^r$ which contains $p$. We say that $K$ is {\it non essentially affine\/} if ${\tilde\psi}^{''}(x)\ne0$ for some $x\in h^{-1}(K)$.

In [Mo], we use the Scale Recurrence Lemma of [MY] in order to prove the following

\vskip .1in

\noindent
{\bf Theorem.} {\it If $K$ and $K'$ are regular Cantor sets of class $C^2$ and $K$ is non essentially affine, then $HD(K+K')=\min\{HD(K)+HD(K'), 1\}$\/}.

This result will be a central tool in the proof of Theorem 1.

\section{Regular Cantor sets defined by the Gauss map}

The Gauss map is the map $g\colon (0,1]\to[0,1]$ given by 
$$g(x)=\{\frac1x\}=\frac1x-\lfloor\frac1x\rfloor, \,\forall\, x\in(0,1].$$
It acts as a shift on continued fractions: if $a_n\in \nb^*, \forall n\ge 1$ then $g([0; a_1, a_2, a_3, \dots]=[0; a_2, a_3, a_4,\dots]$.

Regular Cantor sets defined by the Gauss map (or iterates of it) restricted to some finite union of intervals are closely related to continued fractions with bounded partial quotients. We will often consider such regular Cantor sets associated to {\it complete shifts\/}. A complete shift is associated to finite sets of finite sequences of positive integers in the following way: given a finite set $B=\{\be_1,\be_2,\dots,\be_m\}$, $m\ge2$, where $\be_j\in(\nb^*)^{r_j}$, $r_j\in\nb^*$, $1\le j\le m$ and $\be_i$ does not begin by $\be_j$ for $i \ne j$, the complete shift associated to $B$ is the set $\Sigma(B)\subset(\nb^*)^{\nb}$ of the infinite sequences obtained by concatenations of elements of $B$: 
$$\Sigma(B)=\{(\al_0,\al_1,\al_2,\dots) \, \mid \, \al_j\in B,\,\, \forall\, j\in \nb\}.$$
 Here (and in the rest of the paper), we use the following notation for concatenations of finite sequences: if $\al_j=(a_j^{(1)},a_j^{(2)},\dots,\al_j^{(m_j)})$ then $(\al_0,\al_1,\al_2,\dots)$ means the sequence $(a_0^{(1)},a_0^{(2)},\dots,\al_0^{(m_0)},a_1^{(1)},a_1^{(2)},\dots,\al_1^{(m_1)},a_2^{(1)},a_2^{(2)},\dots,\al_2^{(m_2)},\dots)$. In some cases, when there is no ambiguity, we will write $\al_0 \al_1 \al_2\cdots$ and also $\al_0^N$ to represent the concatenation of $N$ copies of $\al_0$. In some cases some of the $\al_j$ are finite sequences and some cases are single numbers, which are viewed as one-element sequences. Associated to $\Sigma(B)$ is the Cantor set $K(B)\subset[0,1]$ of the real numbers whose continued fractions are of the form $[0;\ga_1,\ga_2,\ga_3,\dots]$, where $\ga_j\in B$, $\forall\, j\ge1$. This is a regular Cantor set. Indeed, if $a_j$ and $b_j$ are respectively the smallest and largest elements of $K(B)$ whose continued fractions begin by $[0; \be_j]$, for $1\le j\le m$, and $I_j=[a_j,b_j]$, then $K(B)$ is the regular Cantor set defined by the map $\psi$ with domain $\bigcup\limits_{j=1}^m I_j$ given by $\psi|_{I_j}=g^{r_j}$, $1\le j\le m$.

We have the following

\vskip .1in

\noindent
{\bf Proposition 1.} {\it The Cantor sets $K(B)$ defined by the Gauss map associated to complete shifts are non essentially affine\/}.

\vskip .1in

\noindent
{\bf Proof:} Let $B=\{\be_1,\be_2,\dots,\be_m\}$, $\be_j=(b_1^{(j)},b_2^{(j)},\dots,b_{r_j}^{(j)})\in(\nb^*)^{r_j}$, $1\le j\le m$. For each $j\le m$, let $x_j=[0; \be_j,\be_j,\be_j,\dots]\in I_j$ be the fixed point of $\psi|_{I_j}=g^{r_j}$. Notice that, since $\be_i$ does not begin by $\be_j$ for $i \ne j$, the $x_j, 1\le j\le m$ are all distinct. Moreover, according to the classical theory of continued fractions, if $p_k^{(j)}/q_k^{(j)}:=[0;b_1^{(j)},b_2^{(j)},\dots,b_k^{(j)}]$, for $1\le j\le m$, $1\le k\le r_j$, we have $I_j\subset\{[0;\be_j,\al],\,\, \al\ge 1\}$ and $\psi|_{I_j}(x)$ is given by
$$
\psi|_{I_j}(x)=\frac{q_{r_j}^{(j)}x-p_{r_j}^{(j)}}{-q_{r_j-1}^{(j)}x+p_{r_j-1}^{(j)}}
$$
(see the appendix); so $x_j$ is the positive root of $q_{r_j-1}^{(j)}x^2+(q_{r_j}^{(j)}-p_{r_j-1}^{(j)})x-p_{r_j}^{(j)}$ (since $x_j$ is the fixed point of $\psi|_{I_j}$).

For each $j\le m$, since $\psi|_{I_j}$ is a M\"obius function with a hyperbolic fixed point $x_j$, there is a M\"obius function $\al_j(x)=\frac{a_jx+b_j}{c_jx+d_j}$ with $\al_j(x_j)=x_j$, $\al_j'(x_j)=1$ such that $\al_j\circ(\psi|_{I_j})\circ \al_j^{-1}$ is an affine map. If we show that the M\"obius functions $\al_1\circ(\psi|_{I_2})\circ \al_1^{-1}$ is not affine then we are done, since the second derivative of a non-affine M\"obius function never vanishes.

Suppose by contradiction that $\al_1\circ(\psi|_{I_2})\circ\al_1^{-1}$ is affine. Since $\al_1\circ(\psi|_{I_1})\circ \al_1^{-1}$ is also affine these two functions have a common fixed point at $\infty$, so $\al_1^{-1}(\infty)=-d_1/c_1$ is a common fixed point of $\psi|_{I_2}$ and $\psi|_{I_1}$, which implies that $\al_1^{-1}(\infty)$ is a common root of $q_{r_1-1}^{(1)}x^2+(q_{r_1}^{(1)}-p_{r_1-1}^{(1)})x-p_{r_1}^{(1)}$ and $q_{r_2-1}^{(2)}x^2+(q_{r_2}^{(2)}-p_{r_2-1}^{(2)})x-p_{r_2}^{(2)}$. Since these polynomials of $\qb[x]$ are irreducible (indeed their roots $x_1$ and $x_2$ are irrational because their continued fractions expansions are infinite), they must be associates in $\qb[x]$, and so their remaining roots $x_1$ and $x_2$ must coincide, which is a contradiction. \qed

\vskip .2in
\noindent{\bf Corollary 1.} $HD(K(B)+K(B'))=\min\{1,HD(K(B))+HD(K(B'))\}$, for every sets $B$, $B'$ of finite sequences of positive integers.

\vskip .2in
\noindent
{\bf Definition 2:} If $\be=(b_1,b_2,\dots,b_{n-1},b_n)$, then $\be^t:=(b_n,b_{n-1},\dots,b_2,b_1)$. Given a set of finite sequences $B$, we define $B^t:=\{\be^t, \be\in B\}$.

\vskip .2in

\noindent
{\bf Proposition 2.} $HD(K(B))=HD(K(B^t))$, for any finite set $B$ of finite sequences.

\vskip .1in

\noindent
{\bf Proof:} This follows from $q_n(\be)=q_n(\be^t)$, \,\, $\forall \be$ (see the appendix of [CF] on properties of continuants), and from the fact that, if $\psi|_{I_j(x)}=\frac{q_{n}^{(j)}x-p_{n}^{(j)}}{-q_{n-1}^{(j)}x+p_{n-1}^{(j)}}$, then 
$$
\psi'|_{I_j(x)}=\frac{-(p_{n}^{(j)}q_{n-1}^{(j)}-p_{n-1}^{(j)}q_{n}^{(j)})}{(-q_{n-1}^{(j)}x+p_{n-1}^{(j)})^2}=\frac{(-1)^n}{(-q_{n-1}^{(j)}x+p_{n-1}^{(j)})^2}
$$
satisfies $(q_{n}^{(j)})^2\le |\psi'|_{I_j(x)}|\le 4(q_{n}^{(j)})^2$, since 
$$\frac1{2q_{n}^{(j)}}\le \frac1{q_{n}^{(j)}+q_{n-1}^{(j)}} \le |q_{n-1}^{(j)}x-p_{n-1}^{(j)}|\le \frac1{q_{n}^{(j)}}.$$
\qed

\vskip .2in

\noindent
{\bf Corollary 2.} $HD(K(B)+K(B^t))=\min\{1,2 \cdot HD(K(B))\}$, for every set $B$ of finite sequences of positive integers.

\section{Fractal dimensions of the spectra}

We recall that the Lagrange spectrum is given by $L=\{\ell(\und\te),\,\, \und\te\in\Sigma\}$, where $\Sigma=(\nb^*)^{\zb}$ and, for $\und\te=(a_n)_{n\in\zb}\in\Sigma$, $\ell(\und\te):=\limsup_{n\to+\infty}(\al_n+\be_n)$, where $\al_n$ and $\be_n$ are defined as the continued fractions $\al_n:=[a_n; a_{n+1},a_{n+2},\dots]$ and $\be_n:=[0;a_{n-1},a_{n-2},\dots]$, while the Markov spectrum is given by $M=\{m(\und\te),\und\te\in\Sigma\}$, where $m(\und\te)=\sup\{\al_n+\be_n,\,\, n\in\zb\}$.

Given a finite sequence $\al=(a_1,a_2,\dots,a_n)\in(\nb^*)^n$, we define its {\it size\/} by $s(\al):=|I(\al)|$, where $I(\al)$ is the interval $\{x\in[0,1] \mid x=[0; a_1,a_2,\dots,a_n,\al_{n+1}]$, $\al_{n+1}\ge 1\}$. If we take $p_0=0$, $q_0=1$, $p_1=1$, $q_1=a_1$ and, for $k\ge 0$, $p_{k+2}=a_{k+2}p_{k+1}+p_k$ and $q_{k+2}=a_{k+2} q_{k+1}+q_k$, then $I(\al)$ is the interval with extremities $[0;a_1,a_2,\dots,a_n]=p_n/q_n$ and $[0;a_1,a_2,\dots,a_{n-1},a_n+1]=\frac{p_n+p_{n-1}}{q_n+q_{n-1}}$ and so
$$
s(\al)=\left| \frac{p_n}{q_n}-\frac{p_n+p_{n-1}}{q_n+q_{n-1}} \right| = \frac1{q_n(q_n+q_{n-1})},
$$
since $p_nq_{n-1}-p_{n-1}q_n=(-1)^{n-1}$. We define $r(\al)=\lfloor\log s(\al)^{-1} \rfloor$ which controls the order of magnitude of the size of $I(\al)$. We also define, for $r \in \nb, P_r= \{\al=(a_1,a_2,\dots,a_n) \mid r(\al) \ge r$ and $r((a_1,a_2,\dots,a_{n-1}))<r\}$.

Write $\Sigma=\Sigma^- \times \Sigma^+$, where $\Sigma^-=(\nb^*)^{\zb_-}$ and $\Sigma^+=(\nb^*)^{\nb}$, and let $\sigma\colon\Sigma\to\Sigma$ be the shift given by $\sigma((a_n)_{n\in\zb})=(a_{n+1})_{n\in\zb}$. We will work with a one-parameter family of subshifts of $\Sigma $ given by $\Sigma_t=\{\und\te\in\Sigma \mid m(\und\te)\le t\}$ for $t\in\re$ (in fact we will take $t\ge3$). Note that $\Sigma_t$ is invariant by transposition and by $\sigma$.

Note that if $\und\te=(a_n)_{n\in\zb}\in\Sigma$ then $\alpha_n+\beta_n>\alpha_n\ge a_n$ for every $n$, and so $m(\und\te)>\sup\{a_n,n\in\zb \}$.  So, if $m(\und\te)\le t$ we have $a_n\le \lfloor t \rfloor$, $\forall\,\, n\in\nb$.
Given $t\in[3,+\infty)$ and $r\in\nb$, let $T:=\lfloor t \rfloor$ and $C(t,r)$ be the set $\{\al=(a_1,a_2,\dots,a_n) \in P_r \mid K_t\cap I(\al)\ne\emptyset\}$. Here $K_t:=\{[0;\ga]| \ga\in\pi_+(\Sigma_t)\}$, where $\pi_+\colon\Sigma\to\Sigma^+$ is the projection associated to the decomposition $\Sigma=\Sigma^-\times\Sigma^+$. Since $\Sigma_t$ is invariant by transposition and by $\sigma$, $K_t$ is invariant by the Gauss map $g$ and $M\cap (-\infty, t) \subset (\nb^* \cap [1, T])+K_t+K_t$. We define $N(t,r):=|C(t,r)|$, where $|\cdot|$ denotes cardinality. Notice that if $r\le s$ then $N(t,r)\le N(t,s)$ and, if $t\le \tilde t$, then $N(t,r)\le N(\tilde t,r)$.

For any finite sequences $\al,\be$ and any positive integers $k_1, k_2\le T$ we have $r(\al\be k_1 k_2)\ge r(\al)+r(\be)$ (see the appendix), so if $C(t,r)=\{\al_1,\al_2,\dots,\al_u\}$ and $C(t,s)=\{\be_1,\be_2,\dots,\be_v\}$, we may cover $K_t$ by
the $T^2uv=T^2N(t,r)N(t,s)$ intervals $I(\al_i\be_j k_1 k_2)$, $1\le i \le u$, $1\le j\le v$, $1\le k_1, k_2\le T$, which satisfy $r(\al_i\be_j k_1 k_2)\ge r+s$, $\forall\,\, i,j,k$. Replacing, if necessary, some of these intervals by larger intervals $I(\ga)$ in $P_{r+s}$, we conclude that $N(t,r+s)\le T^2N(t,r)N(t,s)$ and so
$$
\log(T^2N(t,r+s))\le \log(T^2N(t,r))+\log(T^2N(t,s)),\quad \forall\,\, r,s.
$$
This implies that
$$
\lim_{m\to\infty}\frac1m \log(T^2N(t,m))=\inf_{m\in\nb^*}\frac1m \log(T^2N(t,m))=\lim_{m\to\infty}\frac1m\log(N(t,m))
$$
exists. We will call this limit $D(t)$ (which coincides with the (upper) box dimension of $K_t$, as follows easily from its definition). Notice that $D(t)$ is a non-decreasing function. We will see in the proof of Theorem 1 that $D(t)$ is continuous and that $HD(k^{-1}(-\infty,t))=D(t)$.

\vskip .2in

\noindent
{\bf Lemma 1.} $D(t)$ is right-continuous: given $t_0\in [3,+\infty)$ and $\eta>0$ there is $\delta>0$ such that for $t_0\le t\le t+\delta$ we have $D(t_0)\le D(t) \le D(t+\delta)<D(t_0)+\eta$.

\vskip .1in

\noindent
{\bf Proof:} If for every $t > t_0$, $r$ large, $\frac{\log N(t,r)}{r} \ge D(t_0)+\eta$ we would have $D(t_0) \ge D(t_0) + \eta$, contradiction
(indeed $C(t,r) \subset C(s,r)$ for $t \le s$, and, by compacity, $C(t_0,r) =
{\bigcap}_{t>t_0}\,C(t,r)$). \qed
\vskip .2in

\noindent
{\bf Lemma 2.} {\it Given $t\in (3,+\infty)$ and $\eta \in (0,1)$ there is $\delta>0$ and a Cantor set $K(B)$ defined by the Gauss map associated to a complete shift $\Sigma(B)\subset\Sigma$ such that $\Sigma(B)\subset\Sigma_{t-\delta}$ and
$HD(K(B))>(1-\eta)D(t)$.\/}

\vskip .2in

\noindent
Since the proof of this Lemma is somewhat technical, we will postpone it to Section 6.
\vskip .2in

\noindent
{\bf Lemma 3.} {\it Given a complete shift $\Sigma(X)\subset \Sigma$ (where X is a finite set of finite sequences of positive integers), we have
$$HD(\ell(\Sigma(X)))=HD(m(\Sigma(X)))=\ov{\dim}(\ell(\Sigma(X)))=\ov{\dim}(m(\Sigma(X)))=\min \{2 \cdot HD(K(X)), 1 \}.$$}

\vskip .1in

\noindent
{\bf Proof:}  Let $T$ be the largest element of a sequence in $X$. First of all we clearly have 
$$\ell(\Sigma(X)) \subset m(\Sigma(X))\subset \bigcup\limits_{1\le a\le T\atop 1\le i, j\le R} (a+g^i(K(X))+g^j(K(X))),\;\;\text{where $R$ is the length of}$$
the largest word of $X$, so $HD(\ell(\Sigma(X)) \le HD(m(\Sigma(X)) \le \min \{2 \cdot HD(K(X)), 1 \}$.

Let $\ve>0$ be given. We will show that there are regular Cantor sets $K, K'$ defined by iterates of the Gauss map with $HD(K),HD(K')>HD(K(X))-\ve$ such that $K+K' \subset \ell(\Sigma(X)) \subset m(\Sigma(X))$. Since, by the dimension formula stated in section 2, $HD(K+K')=\min\{HD(K)+HD(K'), 1\}>\min \{2 \cdot HD(K(X)), 1 \}-2\ve$ and $\ve>0$ is arbitrary, the result will follow.

Given a positive integer $n$, let $X^n=\{(\ga_1,\ga_2,\dots,\ga_n)|\ga_j \in X, \forall j \le n \}$. We have $\Sigma(X^n)=\Sigma(X)$ and $K(X^n)=K(X)$. Replacing $X$ by $X^n$ for some $n$ large, we may assume without loss of generality that for any $A \subset X$ (resp. $A^t \subset X^t$) with $|A|\le 2$ (resp. $|A^t|\le 2$), we have $HD(K(X\setminus A))>HD(K(X))-\ve$ (resp. $HD(K(X^t\setminus A^t))>HD(K(X^t))-\ve=HD(K(X))-\ve$).

Order $X$ and $X^t$ in the following way: given $\ga, \tilde \ga \in X \, (\text{resp. }\ga, \tilde \ga\in X^t)$, we say that $\ga < \tilde \ga$ if and only if $[0;\ga]<[0;\tilde \ga]$.

Suppose that the maximum of $m(\Sigma(X))$
is attained at $\tilde \ovt = (\dots,\tilde\ga_{-1},\tilde\ga_0,\tilde\ga_1,\dots),\tilde\ga_i \in X, \forall i \in \zb$, in a position belonging to the sequence $\tilde\ga_0$. Let $X^* = X\backslash\{\min X, \max X\}$,\,\,\, $(X^t)^* = X^t\backslash\{\min X^t, \max X^t\}$. Essentially, $K(X^*)$ and $K((X^t)^*)$ will be the required Cantor sets, but first we have to control the positions where the $\limsup$ is attained (the idea is somewhat similar to the proof that Hall's theorem ([H]) on sums of continued fractions with coefficients bounded by $4$ implies that the Lagrange spectrum contains $[6,+\infty)$) and which words can appear in the beginning of the elements.

For each positive integer $m$, let $C^m$ be the set of sequences $$(\dots, \ga_{-m-2}, \ga_{-m-1}, \tilde\ga_{-m}, \tilde\ga_{-m+1}, \dots, \tilde\ga_{-1}, \tilde\ga_0, \tilde\ga_1, \dots, \tilde\ga_{m-1}, \tilde\ga_m, \ga_{m+1}, \ga_{m+2}, \dots)$$ where $\ga_k \in X^*$ for $k \ge m+1$, \,\,\, $\ga_k^t \in (X^t)^*$ for $k \le -m-1$. Then, for $m$ large enough, there is $\eta>0$ such that for each $\ovt \in C^m$, $\text{sup}(\al_n+\be_n) =m(\ovt)$ is attained only for values of $n$ corresponding to the piece $\tau =(\tilde\ga_{-m}, \tilde\ga_{-m+1}, \dots, \tilde\ga_{-1}, \tilde\ga_0, \tilde\ga_1, \dots, \tilde\ga_{m-1}, \tilde \ga_m)$ of $\ovt$ and, if $n$ does not correspond to the piece $\tau$, then  $\al_n+\be_n<m(\ovt)-\eta$. Indeed, if it is not the case, we may assume without loss of generality that there is a sequence $(m_k)$ tending to $+\infty$ and, for each $k$, $\ovt^{(k)} \in C^{m_k}$ and $n_k$ corresponding to a piece $\ga_{r(k)}$, with $r(k)>m_k$ such that $\al_{n_k}(\ovt^{(k)})+\be_{n_k}(\ovt^{(k)})>m(\ovt^{(k)})-1/k$. Since $\ovt^{(k)}$ converges to $\tilde \ovt$, $m(\ovt^{(k)})$ converges to $m(\tilde \ovt)$ and, by compacity,  if $N_k$ denotes the size of the sequence $\tilde\ga_0, \tilde\ga_1, \dots, \tilde\ga_{m_k-1}, \tilde\ga_{m_k}, \ga_{m_k+1}, \ga_{m_k+2}, \dots, \ga_{r(k)-1}$, $(\sigma^{N_k}(\ovt^{(k)}))$ has a subsequence which converges to some $\hat \ovt = (\dots,\hat\ga_{-1},\hat\ga_0,\hat\ga_1,\dots) \in \Sigma(X)$, with
\,\,\,$\hat\ga_i \in X^*, \forall i \ge 0$, such that $\text{sup}(\al_n+\be_n)=m(\hat \ovt)=m(\tilde \ovt)$ is attained for some $n$ corresponding to the piece $\hat \ga_0$. This is a contradiction, since $m(\tilde \ovt)$ is the maximum of $m(\Sigma(X))$ and, changing $\hat \ga_1$ by $\min X$ or $\max X$, we strictly increase the value of $m(\hat \ovt)$. Notice that the same argument shows that for any $\ovt \in C^m$ and $\ovt^* \in \Sigma(X^*)$, we have $m(\ovt^*)<m(\ovt)-\eta$ (for $m$ large enough).

Now, fixing $m$ with the above properties and $\ga^{(0)} \in X$ such that $(\ga^{(0)})^t \in (X^t)^*$, we may associate to each $x=[0;\ga_1(x),\ga_2(x),\ga_3(x),\dots] \in K(X^*)$ an element $\underline{\Theta}(x) \in C^m$ given by $$\underline{\Theta}(x)=(\dots, \ga^{(0)}, \ga^{(0)}, \tilde\ga_{-m}, \tilde\ga_{-m+1}, \dots, \tilde\ga_{-1}, \tilde\ga_0, \tilde\ga_1, \dots, \tilde\ga_{m-1}, \tilde\ga_m, \ga_1(x), \ga_2(x), \dots)=$$
$$=(\dots, \ga^{(0)}, \ga^{(0)}, \tau, \ga_1(x), \ga_2(x), \dots).$$
For each position $n$ corresponding to the piece $\tau$ of $\underline{\Theta}(x)$, we write $g_n(x)=\al_n(\underline{\Theta}(x))+\be_n(\underline{\Theta}(x))$; in fact $\be_n(\underline{\Theta}(x))$ does not depend on $x$, so, for distinct values of $n$, the functions $g_n$ are distinct rational maps of $x$. This implies that, except for finitely many values of $x$, the values of $g_n(x)$ for these values of $n$ are all distinct. Let $x^{\#}=[0;\ga_1^{\#},\ga_2^{\#},\ga_3^{\#},\dots]$ be one of these values. Since $\text{sup}(\al_n+\be_n) =m(\underline{\Theta}(x^{\#}))$ is attained for values of $n$ corresponding to the piece $\tau$ of $\underline{\Theta}(x^{\#})$, let $n_0$ be the position in $\tau$ for which $m(\underline{\Theta}(x^{\#}))=\al_{n_0}(\underline{\Theta}(x^{\#}))+\be_{n_0}(\underline{\Theta}(x^{\#}))$. For $N$ large enough, taking $\tau^{\#}=((\ga^{(0)})^N, \tau, \ga_1^{\#},\ga_2^{\#},\dots, \ga_N^{\#})$, the following holds: if $\ovt=(\dots, \ga_{-2}, \ga_{-1}, \tau^{\#}, \ga_{1}, \ga_{2}, \dots)$, with $\ga_k \in X^*, (\ga_{-k})^t \in (X^t)^*, \forall k \ge 1$, writing $\tau^{\#}=({\overline a}_{-N_1},\dots,{\overline a}_{-2}, {\overline a}_{-1}, {\overline a}_0, {\overline a}_1, {\overline a}_2, \dots, {\overline a}_{N_2})$, where ${\overline a}_0$ is in the position $n_0$ of $\tau$, we have $m(\ovt)=[{\overline a}_0; {\overline a}_1, {\overline a}_2, \dots, {\overline a}_{N_2},\ga_1, \ga_2, \ga_3, \dots]+[0; {\overline a}_{-1}, {\overline a}_{-2}, \dots, {\overline a}_{-N_1}, {\ga_{-1}}^t, {\ga_{-2}}^t, {\ga_{-3}}^t, \dots]$. It follows that, defining
$$K:=\{[{\overline a}_0; {\overline a}_1, {\overline a}_2, \dots, {\overline a}_{N_2},\ga_1, \ga_2, \ga_3, \dots]| \ga_j \in X^*, \forall j \ge 1\} \text{      and}$$
$$K':=\{[0; {\overline a}_{-1}, {\overline a}_{-2}, \dots, {\overline a}_{-N_1}, {\ga'_1}^t, {\ga'_2}^t, {\ga'_3}^t, \dots]| {\ga'_j}^t \in (X^t)^*, \forall j \ge 1\},$$
we have  $K+K' \subset \ell(\Sigma(X))$. In order to show this, given \\ $x=[{\overline a}_0; {\overline a}_1, {\overline a}_2, \dots, {\overline a}_{N_2},\ga_1, \ga_2, \ga_3, \dots] \in K$ and \\ $y=[0; {\overline a}_{-1}, {\overline a}_{-2}, \dots, {\overline a}_{-N_1}, {\ga'_1}^t, {\ga'_2}^t, {\ga'_3}^t, \dots]\in K'$,\\ and defining, for each positive integer $m$, $\tau^{(m)}=(\ga'_m,\ga'_{m-1},\dots,\ga'_1, \tau^{\#}, \ga_1,\ga_2,\ldots,\ga_m)$, we have, for  
$$\underline{\Theta}^*(x,y)= (\dots, \ga^{(0)}, \ga^{(0)}, \tau^{(1)}, \tau^{(2)}, \tau^{(3)}, \dots),\text{ and }\hat{\underline{\Theta}}(x,y)=(\dots, \ga'_3, \ga'_2, \ga'_1, \tau^{\#}, \ga_1, \ga_2, \ga_3, \dots),$$
$$\ell(\underline{\Theta}^*(x,y))=m(\hat{\underline{\Theta}}(x,y))=x+y.$$
Indeed, there is a sequence of positions $(s_k)$ with  $s_k$ corresponding to the piece $\tau^{(k)}$ of $\underline{\Theta}^*(x,y)$ such that $\sigma^{s_k}(\underline{\Theta}^*(x,y))$ converges to  $\sigma^{n_0}(\hat{\underline{\Theta}}(x,y))$, so $\al_{s_k}(\underline{\Theta}^*(x,y))+\be_{s_k}(\underline{\Theta}^*(x,y))$ converges to $\al_{n_0}(\hat{\underline{\Theta}}(x,y))+\be_{n_0}(\hat{\underline{\Theta}}(x,y))=m(\hat{\underline{\Theta}}(x,y))=x+y$, and, in particular, $\ell(\underline{\Theta}^*(x,y)) \ge m(\hat{\underline{\Theta}}(x,y))=x+y$. On the other hand, there are increasing sequences $(m_k)$ and $(r_k)$ such that the position $m_k$ corresponds to the piece $\tau^{(r_k)}$ in $\underline{\Theta}^*(x,y)$ and $\al_{m_k}(\underline{\Theta}^*(x,y))+\be_{m_k}(\underline{\Theta}^*(x,y))$ converges to  $\ell(\underline{\Theta}^*(x,y))$. Now, if $|m_k-s_{r_k}|$ has a bounded subsequence, then there is $b \in \zb$ such that $\sigma^{m_k}(\underline{\Theta}^*(x,y))$ has a subsequence converging to $\sigma^b(\hat{\underline{\Theta}}(x,y))$, so $\ell(\underline{\Theta}^*(x,y)) =\lim (\al_{m_k}(\underline{\Theta}^*(x,y))+\be_{m_k}(\underline{\Theta}^*(x,y))) \le  m(\hat{\underline{\Theta}}(x,y))=x+y$. On the other hand, if $|m_k-s_{r_k}|$ is unbounded, there is $c \in \zb$ and a subsequence of $\sigma^{m_k}(\underline{\Theta}^*(x,y))$ which converges to $\sigma^c(\ovt^*)$, where $\ovt^*$ is an element of $\Sigma(X^*)$, but in this case we would have $\ell(\underline{\Theta}^*(x,y)) \le m(\ovt^*)<m(\hat{\underline{\Theta}}(x,y))-\eta$, which is a contradiction.

Finally, notice that $K$ and $K'$ are diffeomorphic respectively to $K(X^*)$ and $K((X^t)^*)$, so $HD(K)=HD(K(X^*))>HD(K(X))-\ve$ and $HD(K')=HD(K((X^t)^*))>$\break $HD(K(X^t))-\ve=HD(K(X))-\ve$. \qed

\section{Proofs of the main results}

\noindent{\bf Proof of Theorem 1:} Applying Lemma 3 to the complete shift $\Sigma(B)$ obtained in Lemma 2, we get that, for any $\eta>0$, there is $\delta>0$ such that \hfill\break $\min\{2(1-\eta)D(t),1\} \le \min\{2HD(K(B)),1\} \le HD(L\cap(-\infty,t-\delta]) \le HD(L \cap (-\infty,t)) \le HD(M \cap (-\infty,t)) \le \ov{\dim}(M\cap(-\infty,t))\le \min\{2 \cdot HD(K_t),1\} \le \min\{2 \cdot D(t),1\}$ \break (since $\ell(\Sigma(B))\subset L\cap(-\infty,t-\delta]$, $L\cap (-\infty, t) \subset M\cap (-\infty, t) \subset (\nb^* \cap [1, T])+K_t+K_t$ and $D(t)$ is the upper box dimension of $K_t$), and so, if
$d(t):=HD(L\cap(-\infty,t))$, we have $d(t)= HD(M\cap(-\infty,t))=\ov{\dim}(L\cap(-\infty,t))=\ov{\dim}(M\cap(-\infty,t))=\min\{2 \cdot D(t),1\}$.

In order to conclude the proof of the first assertion of i), it is enough to show that $HD(k^{-1}(-\infty,t))=D(t)$. In the notation of Lemma 3, let $x\in K, y\in K'$. For each $z=[0; \al_1, \al_2, \dots]\in K(X^*)$, define
\begin{multline*}
\und\la(z)={\und\la}_{x,y}(z)=(\al_{1!},\tau^{(1)},\al_{2!},\tau^{(2)},\al_3,\al_4,\al_5,\al_{3!},\tau^{(3)},\al_7,\dots,\al_{4!},\tau^{(4)},\\
\al_{25},\al_{26},\dots,\al_{5!},\tau^{(5)},\;\dots\;,\al_{r!},\tau^{(r)},\al_{r!+1},\dots),
\end{multline*}
and $h(z)=[0;\und \la(z)]$. We have, as before, $k(h(z))=x+y$. On the other hand, given any $\rho>0$, we have $|z-z'|=O(|h(z)-h(z')|^{1-\rho})$ for $|z-z'|$ small, so $HD(k^{-1}(x+y)) \ge HD(K(B^*))>HD(K(B))-\ve$. As before, we get 
$HD(k^{-1}(-\infty,t)) \ge HD(k^{-1}(-\infty,t-\delta])$ $\ge HD(k^{-1}(x+y)) > HD(K(B))-\ve > (1-\eta)D(t)-\ve.$ 
Since $\eta$ and $\ve$ are arbitrary, $HD(k^{-1}(-\infty,t)) \ge D(t)$. For the reverse inequality, let $w \in  k^{-1}(-\infty,t)$. We have $\limsup_{n\to\infty}(\al_n(w)+\be_n(w))=k(w)<t$, so there is $n_0 \in \nb$ such that $n \ge n_0 \implies \al_n(w)+\be_n(w)<t$. This implies that $k^{-1}(-\infty,t) \subset {\bigcup}_{n \in \nb}(g^{-n}(K_t))$, where $g$ is the Gauss map, so $HD(k^{-1}(-\infty,t)) \le D(t)$. Thus we have $HD(k^{-1}(-\infty,t))=D(t)$. 

Recall that, by Lemma 1, $D(t)$ is a right-continuous function.

Thus we have
\begin{eqnarray*}
D(t)=HD(k^{-1}(-\infty,t)) &\le& HD(k^{-1}(-\infty,t]) \le \lim\limits_{s\to t+} HD(k^{-1}(-\infty,s)) \\
&=& \lim\limits_{s\to t+} D(s) = D(t).
\end{eqnarray*}
Then $HD(k^{-1}(-\infty,t])=HD(k^{-1}(-\infty,t))=D(t)$, and we conclude that $d(t) =\break \min\{2HD(k^{-1}(-\infty,t)),1\} = \min\{2HD(k^{-1}(-\infty,t]),1\}$.

Finally, $D(t)$ is left-continous (and so is continuous), since, by Lemma 2, given $t\in [3,+\infty)$ and $\eta\in (0,1)$, there is $\delta>0$ such that $D(t-\delta)\ge HD(K(B))>(1-\eta)D(t)$, so $\lim\limits_{s\to t-} D(s) = D(t)$.

In order to conclude, notice that, for each positive integer $m$, \break $\Sigma(\{21^{2m}2,21^{2m+2}2\}) \subset \Sigma_{3+2^{-m}}$ (notice that $[2;1, 1,  1,\ldots]+[0;2, 1, 1, 1,\ldots]=3$), so \break $D(3+\ve)>0$ for every $\ve>0$ and, since $\Sigma_{\sqrt {12}}=\{1,2\}^{\zb}$, $D(\sqrt{12})=HD(K_{\sqrt{12}})=HD(K(\{1,2\}))=HD(C_2)=0,53128\ldots>1/2$, so, since $D(t)$ is continuous, there is $\delta>0$ such that $D(\sqrt{12}-\delta)>1/2$, and thus we have $d(\sqrt {12}-\delta)=\min\{2 \cdot D(\sqrt {12}-\delta),1\}=1$.
\qquad \qed

\noindent{\bf{Remark:}}
It follows from the above proof and from the general estimates of fractal dimensions of regular Cantor sets of Chapter 4 of [PT] that there is a constant $C>0$ such that, for each positive integer $m$, $D(3+2^{-m})\ge HD(K(\{21^{2m}2,21^{2m+2}2\}))>C/m$. This gives another proof of the fact that the functions $D(t)$ and $d(t)$ are not H\"older continuous.
\vglue .2in

\noindent{\bf Proof of Theorem 2:}
Given $m \ge 2$, let $C_m=\{\al=[0;a_1,a_2,a_3,\ldots] \in [0,1]| a_k \le m, \forall k \ge 1\}$. M. Hall proved in [H] that $C_4+C_4=\{\al+\be| \al,\be \in C_4\}=[\sqrt{2}-1,4(\sqrt{2}-1)]$. On the other hand, we have $\lim_{m\to\infty}HD(C_m)=1$. In fact, Jarn\'\i k proved in [J] that 
$$HD(C_m)>1-\frac{1}{m \cdot \log 2}, \forall m>8.$$

Let now $t \ge 7$ be given. Let $m=\lfloor t \rfloor -3$. There are an integer $n \in \{m+2,m+3\}$ and $\al=[0;a_1,a_2,a_3,\ldots], \be=[0;b_1,b_2,b_3,\ldots] \in C_4$ such that $t=n+\al+\be$. For each $r\ge 1$, let $\tilde\tau^{(r)}$ and $\hat\tau^{(r)}$ be respectively the sequences \hfill\break $(m+1,b_{2r-1},b_{2r-2},\ldots,b_2,b_1,n,a_1,a_2,\ldots,a_{2r-2},a_{2r-1},m+1)$ and \hfill\break $(m+1,b_{2r},b_{2r-1},\ldots,b_2,b_1,n,a_1,a_2,\ldots,a_{2r-1},a_{2r},m+1)$. Consider now the maps \break $\tilde h, \hat h: C_m \rightarrow [0,1]$ given by
\begin{multline*}
\tilde h(z) = \tilde h([0;c_1,c_2,c_3,\ldots] )=[0;c_{1!},\tilde\tau^{(1)},c_{2!},\tilde\tau^{(2)}, c_3,c_4,c_5,c_{3!},\tilde\tau^{(3)}, \\
c_7, c_8,\dots,c_{4!},\tilde\tau^{(4)}, c_{25}, \dots,c_{5!}, \tilde\tau^{(5)},\: \dots\:, c_{r!},\tilde\tau^{(r)},c_{r!+1},\dots].
\end{multline*}
\begin{multline*}
\hat h(z) = \hat h([0;c_1,c_2,c_3,\ldots] )=[0;c_{1!},\hat\tau^{(1)},c_{2!},\hat\tau^{(2)}, c_3,c_4,c_5,c_{3!},\hat\tau^{(3)}, \\
c_7, c_8,\dots,c_{4!},\hat\tau^{(4)}, c_{25}, \dots,c_{5!}, \hat\tau^{(5)},\: \dots\:, c_{r!},\hat\tau^{(r)},c_{r!+1},\dots].
\end{multline*}

It is easy to see that $k(\tilde h(z))=t$ for every $z \in C_m$. Moreover, since $[x_0;x_1,x_2,x_3,x_4,\dots]$ is increasing in $x_0, x_2, x_4,\ldots$ and decreasing in $x_1, x_3, x_5, \ldots$, $\tilde h(z)\in \text{Exact}(t^{-1})$ and $\hat h(z)\in \text{Exact}'(t^{-1})$. On the other hand,
given any $\rho>0$, we have $|z-z'|=O(|\tilde h(z)-\tilde h(z')|^{1-\rho})$ and $|z-z'|=O(|\hat h(z)-\hat h(z')|^{1-\rho})$ for $|z-z'|$ small, so $HD(\text{Exact}(t^{-1})) \ge HD(C_m)$ and $HD(\text{Exact}'(t^{-1})) \ge HD(C_m)$. Since $\lim_{m\to\infty}HD(C_m)=1$, we are done.	   \qquad\qed
\\

\noindent{\bf Proof of Theorem 3:}
Let $x \in L'$. Consider a sequence $x_n$ converging to $x$, $x_n \in L$, $x_n \ne x$. Choose $\und{\te}^{(n)} \in \Sigma$
such that $x_n = \ell(\und{\te}^{(n)})$. Let $\und{\te}^{(n)}=(b_j^{(n)})_{j\in\zb}$ and assume $b_j^{(n)} \le 4$, $\forall\, j$, $\forall\, n$ (which is possible since me may assume that the $x_n$ are not in Hall's ray). We have $x_n ={\limsup}_{j\to\infty}(\al_j^{(n)} + \be_j^{(n)})$. Given $\delta >0$, \,\,$\exists\, n_0 \in \nb$ large such that $n \ge n_0 \Rightarrow
|\ell(\und{\te}^{(n)})-x| < \delta$ and there are infinitely many $j \in \nb$ such that
$|\al_j^{(n)} + \be_j^{(n)}-x| < \delta$. Let $N =
\lceil \delta^{-1} \rceil$. Given such a pair $(j,n)$ consider the finite sequence with $2N+1$ terms $(b_{j-N}^{(n)},b_{j-N+1}^{(n)},\dots, b_j^{(n)},\dots, b_{j+N}^{(n)}) =: S(j,n)$. There is a
sequence $S$ such that for infinitely many values of $n$, \,$S$ appears
infinitely many times as $S(j,n), \, j \in \nb$, i.e., there are $j_1(n) < j_2(n) <\dots$
with $\lim_{i\to \infty} (j_{i+1}(n)-j_i(n))=\infty$ and $S(j_i(n),n) = S$, $\forall\,i \ge 1$, for all $n$ in some infinite set $A \subset \nb$.

Consider the sequences $\be(i,n)$ for $i\ge 1$, $n \in A$ given by
$$
\be(i,n) = (b_{j_i(n)+N+1}^{(n)}, b_{j_i(n)+N+2}^{(n)},\dots,
b_{j_{i+1}(n)+N}^{(n)}).
$$
There are $(i_1,n_1)$ and $(i_2,n_2)$ for which there is no sequence $\ga$ such
that $\be(i_1,n_1)$ and $\be(i_2,n_2)$ are concatenations of copies of $\ga$,
otherwise $x_n$ would be constant for $n \in A$. This implies that, taking $B = \{\be(i_1,n_1) \be(i_2,n_2), \be(i_2,n_2) \be(i_1,n_1)\}$, \,\, $K(B)$ is a regular Cantor set, so, as in Lemma 3,
$\ell(K(B))$ contains a regular Cantor set $\hat K$ with $d(x,\hat K) \le 2\delta$.
\qquad\qed

\vglue .2in
\section{Proof of Lemma 2}
Let $\tau=\eta/40$. Since $t>3$, we have $D(t)>0$, and so we may choose $r_0\in\nb$ large such that, for $r\ge r_0$, $|\frac{\log N(t,r)}r-D(t)|<\frac\tau 2 D(t)$. Let $B_0:=C(t,r_0)$ and $N_0:=N(t,r_0)=|B_0|$. Let $k=8N_0^2\lceil 2/\tau\rceil$.
Take $\tilde B=\{\be=\be_1\be_2\cdots \be_k \mid \be_j\in B_0,1\le j\le k $ and $K_t\cap I(\be)\ne \emptyset\}$.

Given $\be=\be_1\be_2\cdots\be_k\in\tilde B$ (with $\be_i\in B_0$, $1\le i\le k$), we say that $j$, $1\le j\le k$, is a {\it right-good\/} position of $\be$ if there are elements $\be^{(s)}=\be_1\be_2\cdots\be_{j-1} \be_j^{(s)} \be_{j+1}^{(s)}\cdots \be_k^{(s)}$, $s=1,2$, of $\tilde B$ such that we have the following inequality of continued fractions: $[0; \be_j^{(1)}]<[0;\be_j]<[0;\be_j^{(2)}]$. We say that $j$ is a {\it left-good\/} position if there are elements $\be^{(s)}=\be_1^{(s)}\be_2^{(s)}\cdots\be_{j-1}^{(s)}\be_j^{(s)}\be_{j+1}\be_{j+2}\cdots\be_k$, $s=3,4$, of $\tilde B$ such that $[0;(\be_j^{(3)})^t]<[0;\be_j^t]<[0;(\be_j^{(4)})^t]$. Finally, we say that $j$ is a {\it good\/} position if it is both right-good and left-good.

We will show that most positions of most words of $\tilde B$ are good. Let us first estimate $|\tilde B|$. It follows from Lemma A2 of the appendix that, for $\be\in\tilde B$, $s(\be)<(2e^{-r_0})^k<e^{-k(r_0-1)}$. Moreover, since $N(t,k(r_0-1))\ge  \frac{1}{T^2} e^{k(r_0-1)D(t)}$, $\{I(\be); \be\in \tilde B\}$ is a covering of $K_t$ by intervals of size smaller than $e^{-k(r_0-1)}$ and the function $h:\tilde B \to C(t,k(r_0-1))$ defined by $h(\be)=h((\be_1\be_2\dots \be_k))=(\be_1\be_2\dots \be_j)$, where $j = \min\{i\; ;i\le k\; \text{ and }\;r((\be_1\be_2\dots \be_i))\ge k(r_0-1)\;\}$ is onto, we have: 

\begin{eqnarray}
|\tilde B|\ge \; \frac{1}{T^2} e^{k(r_0-1)D(t)}&>& 2\,e^{k(r_0-2)D(t)},\quad \text{since $k$ is large} \nonumber\\
&\ge& 2\,e^{(1-\tau/2)r_0 kD(t)},\quad \text{since $r_0$ is large} \nonumber\\
&>& 2\,e^{(1-\tau)(1+\tau/2)r_0kD(t)}\nonumber\\
&>& 2\,N_0^{(1-\tau)k},\quad \text{since}\: N(t,r_0)<e^{\left(1+\frac{\tau}{2}\right)D(t) r_0} .\nonumber\\
\nonumber
\end{eqnarray}

Now, let us estimate the number of words $\be\in\tilde B$ such that at least $k/20$ positions of $\be$ are not right good: we have at most $2^k$ choices for the set of the $m\ge k/20$ positions which are not right-good. Once we choose this set of positions, if $j$ is such a position and $\be_1,\be_2,\dots,\be_{j-1}\in B_0$ are already chosen, there are at most two (the largest and the smallest) choices for $\be_j\in B_0$ such that for some $\be=\be_1\be_2\cdots\be_{j-1}\be_j\be_{j+1}\cdots\be_k\in\tilde B$ the position $j$ is not right good. If $j$ is any other position, we have of course at most $N_0=|B_0|$ possible choices for $\be_j$, so we have at most $2^m\cdot N_0^{k-m}\le 2^{k/20}N_0^{19k/20}$ words in $\tilde B$ with this chosen set of $m$ positions which are not right-good. Therefore, the number of words $\be\in\tilde B$ for which the number of positions which are not right-good is at least $k/20$ is bounded by $2^k\cdot 2^{k/20}\cdot N_0^{19k/20}=2^{21k/20}\cdot N_0^{19k/20}$. Analogously, the number of words $\be\in\tilde B$ for which there are at least $k/20$ positions which are not left-good is also bounded by $2^{21k/20}\cdot N_0^{19k/20}$.

This implies that for at least $|\tilde B|-2\cdot 2^{21k/20}\cdot N_0^{19k/20} > 2N_0^{(1-\tau)k}-2^{1+21k/20}\cdot N_0^{19k/20} > N_0^{(1-\tau)k}$ words of $\tilde B$, the number of good positions is at least $9k/10$.
Let us call such an element of $\tilde B$ an {\it excellent\/} word.

If $\be=\be_1\be_2\cdots\be_k\in\tilde B$ (with $\be_j\in B_0$, $1\le j\le k$) is an excellent word, we may find $\lceil 2k/5\rceil$ positions $i_1,i_2,\dots,i_{\lceil 2k/5\rceil}\le k$ with $i_{s+1}\ge i_s+2$, $\forall\, s< \lceil 2k/5\rceil$, such that the positions $i_1,i_1+1$, $i_2, i_{2}+1,\dots,i_{\lceil 2k/5\rceil}$, $i_{\lceil 2k/5\rceil}+1$ are good. Since $k=8N_0^2\lceil 2/\tau\rceil$, we may take, for $1\le s\le 3N_0^2$, $j_s:=i_{s\lceil 2/\tau\rceil}$ (notice that $3N_0^2\lceil 2/\tau\rceil<\frac{16}5N_0^2\lceil 2/\tau\rceil=2k/5$), so we have $j_{s+1}-j_s\ge 2\lceil 2/\tau\rceil$, $\forall\, s<3N_0^2$, and the positions $j_s,j_s+1$ are good for $1\le s\le 3N_0^2$.

Now, the number of possible choices of $(j_1,j_2,\dots,j_{3N_0^2})$ is bounded by $\binom{k}{3N_0^2}<2^k$ and, given $(j_1,j_2,\dots,j_{3N_0^2})$ the number of choices of $(\be_{j_1},\be_{j_1+1},\dots,\be_{j_{3N_0^2}},\be_{j_{3N_0^2}+1})$ is bounded by $N_0^{6N_0^2}$. So, we may choose $\hat {\j}_1,\hat {\j}_2,\dots,\hat {\j}_{3N_0^2}$ with $\hat {\j}_{s+1}-\hat {\j}_s\ge 2\lceil 2/\tau\rceil$, $\forall\, s<3N_0^2$, and words
$\hat\be_{\hat {\j}_1},\hat\be_{\hat {\j}_1+1}, \hat\be_{\hat {\j}_2},\hat\be_{\hat {\j}_2+1},\dots,\hat\be_{\hat {\j}_{3N_0^2}},\hat\be_{\hat {\j}_{3N_0^2}+1} \in B_0$ such that the set
$X:=\{\be=\be_1\be_2\cdots\be_k\in\tilde B$ excellent$\,|\,\, \hat {\j}_s,\hat {\j}_s+1$ are good positions and $\be_{\hat {\j}_s}=\hat\be_{\hat {\j}_s}, \be_{\hat {\j}_s+1}=\hat\be_{\hat {\j}_s+1}, \forall\, s\le 3N_0^2\}$ has at least $\frac{N_0^{(1-\tau)k}}{2^k \cdot N_0^{6N_0^2}} > N_0^{(1-2\tau)k}$ elements, as $N_0$ and $k$ are large.

Since $N_0=|B_0|$, there are $N_0^2$ possible choices for the pairs $(\hat\be_{\hat {\j}_s},\hat\be_{\hat {\j}_s+1})$. We will consider, for $1\le s<t\le 3N_0^2$, the projections $\pi_{s,t}\colon X\to B_0^{\hat{\j}_t-\hat {\j}_s}$ given by $\pi_{s,t}(\be_1\be_2\cdots\be_k)=(\be_{\hat {\j}_s+1},\be_{\hat {\j}_s+2},\dots,\be_{\hat {\j}_t})$. We will show that the images of many of these projections are large.

For each pair $(s,t)$ with $1\le s<t\le 3N_0^2$ such that $|\pi_{s,t}(X)| < N_0^{(1-10\tau)(\hat {\j}_t-\hat {\j}_s)}$, we will exclude from $\{1,2,\dots,3N_0^2\}$ the indices $s,s+1,\dots,t-1$.
Let us estimate the total number of indices excluded: the set of excluded indices is the union of the intervals $[s,t)$ (intersected with $\zb$) over the pairs $(s,t)$ as above. Now we use the elementary fact that, given a finite family of intervals, there is a subfamily of disjoint intervals whose sum of lenghts is at least half of the measure of the union of the intervals of the original family. We apply this fact to the above intervals $[s,t)$. Suppose that the total number of indices excluded is at least $2N_0^2$. By the above fact, we may find a disjoint collection of intervals $[s,t)$ as above whose sum of lenghts is at least $N_0^2$. Let us call $\cal P$ the set of these pairs $(s,t)$. Since $\hat {\j}_t-\hat {\j}_s\ge 2(t-s)\lceil 2/\tau\rceil, \forall t>s$, the sum of $(\hat {\j}_t-\hat {\j}_s)$ for $(s,t) \in \cal P$ is at least $2N_0^2\lceil 2/\tau\rceil$. Since for each  pair $(s,t) \in \cal P$ we have $|\pi_{s,t}(X)|<N_0^{(1-10\tau)(\hat {\j}_t-\hat {\j}_s)}$, we get
\begin{eqnarray}
N_0^{(1-2\tau)k} <|X| &<& N_0^{(1-10\tau)\sum\limits_{(s,t)\in \cal P}(\hat {\j}_t-\hat {\j}_s)} \cdot N_0^{\#\{i; i\: \notin [\hat {\j}_s,\hat {\j}_t), \forall (s,t) \in \cal P\}} \nonumber\\
&<&N_0^{(1-10\tau)\cdot 2N_0^2\lceil 2/\tau\rceil} \cdot N_0^{k-2N_0^2\lceil 2/\tau\rceil}, \nonumber
\nonumber
\end{eqnarray}
since we have at most $N_0$ choices for $\be_i$ for each index $i$ which does not belong to the union of the intervals $[\hat {\j}_s,\hat {\j}_t)$ associated to these pairs $(s,t)$. However, this is a contradiction, since this inequality is equivalent to $N_0^{20\tau N_0^2\lceil 2/\tau\rceil}<N_0^{2\tau k}$, which cannot hold, because $2\tau k=16\tau N_0^2\lceil 2/\tau\rceil<20\tau N_0^2\lceil 2/\tau\rceil$. So, the total number of excluded indices is smaller than $2N_0^2$.

Now, there are at least $N_0^2+1$ indices which are not excluded. We will have two non-excluded indices $s<t$ such that $\hat\be_{\hat {\j}_s}=\hat\be_{\hat {\j}_t}$ and $\hat\be_{\hat {\j}_s+1}=\hat\be_{\hat {\j}_t+1}$. We claim that, for $B:=\pi_{s,t}(X)$, the shift $\Sigma(B)$ satisfies the conclusions of the statement.

Indeed, since $s$ and $t$ are not excluded, we have $|B|\ge N_0^{(1-10\tau)(\hat {\j}_t-\hat {\j}_s)}$. Moreover, by Proposition A1 of the appendix, for every $\al\in B$ we have
$$|I(\al)|=s(\al)>(2(T+1)^2 e^{r_0})^{-(\hat {\j}_t-\hat {\j}_s)} > e^{-(\hat {\j}_t-\hat {\j}_s)(r_0+\lceil \log(2(T+1)^2) \rceil)}.$$
 So, the Hausdorff dimension of $K(B)$ is at least
$$
\frac{(1-10\tau)\log N_0}{r_0+\lceil \log(2(T+1)^2) \rceil} > \frac{(1-10\tau)r_0}{r_0+\lceil \log(2(T+1)^2) \rceil}\cdot (1-\frac\tau 2)D(t) > \left(1-12\tau\right)D(t)> (1-\eta)D(t).
$$
On the other hand, if $\tilde k:=\hat {\j}_t-\hat {\j}_s$, $\ga_1:=\hat\be_{\hat {\j}_s+1}=\hat\be_{\hat {\j}_t+1}$ and $\ga_2:=\hat\be_{\hat {\j}_t}=\hat\be_{\hat {\j}_s}$, all elements of $B$ are of the form $\ga_1\be_2\be_3\cdots\be_{\tilde k-1}\ga_2$, where $\ga_1,\be_2,\be_3,\dots,\be_{\tilde k-1}$, $\ga_2\in B_0$ and there are $\ga_1',\ga_1''$, $\ga_2',\ga_2''\in B_0$ with $[0;\ga_2']<[0;\ga_2]< [0;\ga_2'']$ and $[0;(\ga_1')^t]<[0;\ga_1^t]<[0;(\ga_1'')^t]$ such that
\begin{eqnarray}
I(\ga_1'\be_2\be_3\cdots\be_{\tilde k-1}\ga_2\ga_1)\cap K_t\ne \emptyset,\quad I(\ga_1''\be_2\be_3\cdots\be_{\tilde k-1}\ga_2\ga_1)\cap K_t\ne \emptyset, \nonumber\\
I(\ga_2\ga_1\be_2\be_3\cdots\be_{\tilde k-1}\ga_2')\cap K_t\ne \emptyset,\quad I(\ga_2\ga_1\be_2\be_3\cdots\be_{\tilde k-1}\ga_2'')\cap K_t\ne \emptyset.\nonumber
\end{eqnarray}
We will show that this implies the existence of $\delta>0$ such that $\Sigma(B)\subset\Sigma_{t-\delta}$. Let $\ga_1^t=(c_1,c_2,\dots,c_{m_1})$, with $c_j\in \nb^*$, $\forall\, j\le m_1$, and $\ga_2=(d_1,d_2,\dots,d_{m_2})$ with $d_j\in\nb^*$, $\forall\, j\le m_2$. Let $\ga_1\be_2\be_3\cdots\be_{\tilde k-1}\ga_2\in B$ where $\be_2\be_3\cdots\be_{\tilde k-1}=a_1a_2\cdots a_{\tilde m}$ with $a_j\in\nb^*$, $\forall\, j\le\tilde m$. We want to estimate three kinds of sums of continued fractions. The first one are sums of continued fractions beginning by $[a_j;a_{j+1},\dots,a_{\tilde m},\ga_2,\ga_1,\dots]+[0;a_{j-1},\dots,a_1,\ga_1^t,\ga_2^t,\dots].$ Let us assume, without loss of generality, that \hfill\break $q_{m_2+\tilde m -j}(a_{j+1},\ldots, a_{\tilde m}, \ga_2) \le q_{m_1+j-1}(a_{j-1},\ldots, a_1, \ga_1^t)$ (the other case, when the reverse inequality $q_{m_1+j-1}(a_{j-1},\ldots, a_1, \ga_1^t) \le q_{m_2+\tilde m -j}(a_{j+1},\ldots, a_{\tilde m}, \ga_2)$ holds, is symmetric).
Assume also that $[a_j; a_{j+1},\ldots,a_{\tilde m},\ga_2]<[a_j; a_{j+1},\ldots,a_{\tilde m},\ga_2']$ (otherwise we change $\ga_2'$ by $\ga_2''$). This allows us to exhibit $\delta > 0$ such that, for any $\und \theta^{(i)} \in \{1, 2, \dots, T\}^{\nb},\; 1 \le i \le 4$,

\begin{eqnarray*}\lefteqn{[a_j;a_{j+1},\dots,a_{\tilde m},\ga_2,\und \theta^{(1)}]+[0;a_{j-1},\dots,a_1,\ga_1^t,\ga_2^t,\und \theta^{(2)}] <}\\
& <[a_j;a_{j+1},\dots,a_{\tilde m},\ga_2',\und \theta^{(3)}]+[0;a_{j-1},\dots,a_1,\ga_1^t,\ga_2^t,\und \theta^{(4)}]-\delta.
\end{eqnarray*}
Indeed, by the Lemma A1 of the appendix,
$$[a_j;a_{j+1},\dots,a_{\tilde m},\ga_2',\und \theta^{(3)}]-[a_j;a_{j+1},\dots,a_{\tilde m},\ga_2,\und \theta^{(1)}]>
\frac1{(T+1)(T+2)q_{m_2+\tilde m -j}(a_{j+1},\ldots, a_{\tilde m}, \ga_2)^2}$$

\begin{multline*}
\text{ and }\,\mid[0;a_{j-1},\dots,a_1,\ga_1^t,\ga_2^t,\und \theta^{(4)}]-[0;a_{j-1},\dots,a_1,\ga_1^t,\ga_2^t,\und \theta^{(2)}]\mid < \\
\frac1{q_{m_1+m_2+j-1}(a_{j-1},\ldots, a_1, \ga_1^t, \ga_2^t)^2}<\frac1{(F_{m_2+1}q_{m_1+j-1}(a_{j-1}, \ldots, a_1, \ga_1^t))^2} \le \\
\frac1{(F_{m_2+1}q_{m_2+\tilde m -j}(a_{j+1},\ldots, a_{\tilde m}, \ga_2))^2}\le \frac1{2(T+1)(T+2)q_{m_2+\tilde m -j}(a_{j+1},\ldots, a_{\tilde m}, \ga_2)^2}\\(\text{here we use that}\; m_2 \text{ is large};\\ (F_n)\text{ denotes Fibonacci's sequence, given by } F_0=0, F_1=1, F_{n+2}=F_{n+1}+F_n, \forall n \ge 0).\\
\end{multline*} 
So, the inequality holds with
$$\delta:=\frac1{(T+1)^{2(m_1+m_2+\tilde m)}}<\frac1{2(T+1)(T+2)q_{m_2+\tilde m -j}(a_{j+1},\ldots, a_{\tilde m}, \ga_2)^2}.$$
On the other hand, $I(\ga_2\ga_1\be_2\be_3\cdots\be_{\tilde k-1}\ga_2')\cap K_t\ne \emptyset$, so there are $\und \theta^{(3)}$ and $\und \theta^{(4)}$ such \break that $(\und \theta^{(4)})^t \ga_2\ga_1 \be_2\be_3\cdots\be_{\tilde k-1}\ga_2' \und\theta^{(3)} \in \Sigma_t$, and thus $[a_j;a_{j+1},\dots,a_{\tilde m},\ga_2',\und \theta^{(3)}]+ \break [0;a_{j-1},\dots,a_1,\ga_1^t,\ga_2^t,\und \theta^{(4)}]\le t$, which implies that, for any $\und \theta^{(i)} \in \{1, 2, \dots, T\}^{\nb}, i=1,2$, \break $[a_j;a_{j+1},\dots,a_{\tilde m},\ga_2,\und \theta^{(1)}]+[0;a_{j-1},\dots,a_1,\ga_1^t,\ga_2^t,\und \theta^{(2)}] < t-\delta$.

The other two kinds of sums of continued fractions we want to estimate are sums of continued fractions beginning by $[d_j; d_{j+1},\ldots,d_{m_2},\ga_1,\ldots]+[0; d_{j-1},\ldots, d_1, a_{\tilde m},\ldots, a_1, \ga_1^t,\ldots]$ and, symmetrically, sums of continued fractions beginning by $[0; c_{j+1},\ldots, c_{m_1}, \ga_2^t,\ldots]+[c_j; c_{j-1},\ldots, c_1, a_1,\ldots, a_{\tilde m}, \ga_2,\ldots]$. We have: $$q_{m_2-j+m_1}(d_{j+1},\ldots,d_{m_2}, \ga_1) \le q_{j-1+\tilde m +m_1}(d_{j-1},\ldots, d_1, a_{\tilde m},\ldots, a_1, \ga_1^t)$$
(indeed, $\tilde m/(m_1+m_2)$ is large when $\eta$ and $\tau$ are small, depending on $T$).
Assume that  $[d_j; d_{j+1},\ldots,d_{m_2}, \ga_1] < [d_j; d_{j+1},\ldots,d_{m_2}, \ga_1']$ (otherwise we change $\ga_1'$ by $\ga_1''$). Since $I(\ga_2\ga_1 \be_2\be_3\cdots\be_{\tilde k-1}\ga_2\ga_1')\cap K_t\ne \emptyset$, estimates analogous to the previous ones \break imply that, for any $\und \theta^{(i)} \in \{1, 2, \dots, T\}^{\nb},\: i=1,2$, we have $[d_j; d_{j+1},\ldots,d_{m_2}, \ga_1,\und \theta^{(1)}]+$\break $[0; d_{j-1},\ldots, d_1, a_{\tilde m},\ldots, a_1, \ga_1^t,\ga_2^t,\und\theta^{(2)}] < t-\delta$.\\
This implies that the complete shift $\Sigma(B)$ satisfies the conditions of the statement, which concludes the proof of the Lemma.
\qed

As we said before, this proof doesn't give any estimate on the modulus of continuity of $d(t)$. Indeed, in the beginning of the proof of Lemma 2, we used the fact that $u(r)=\log(T^2N(t,r))$ is subadditive in order to guarantee the existence of $r_0\in\nb$ large such that, for $r\ge r_0$, $|\frac{\log N(t,r)}r-D(t)|<\frac\tau 2 D(t)$ (recall that $D(t)=\lim_{m\to\infty}\frac1m \log(T^2N(t,m))=\lim_{m\to\infty}\frac1m\log(N(t,m))$). However, this gives no estimate on $r_0$. Consider, for instance, the function $v(n)$ given by $v(n)=2n$ for $n\le M_0$ and $v(n)=n+M_0$ for $n>M_0$, where $M_0$ is a large positive integer. It is subadditive, increasing, $\lim_{n\to \infty} v(n)/n=1$ but $v(M_0)/M_0=2$, and $M_0$ can be taken arbitrarily large. However it is possible to adapt the proof in order to give an estimate on the modulus of continuity of $d(t)$, using an idea of [FMM]. 

Given $\ve>0$ (which we may assume to be smaller than $\frac17<\frac1{10\log2}$), we want to obtain $\delta\in (0,1)$ as an explicit function of $\ve$ such that $D(t-\delta)>D(t)-\ve$. Of course there is no loss of generality in assuming $D(t)\ge \ve$. We may also assume that $T=\lfloor t \rfloor < 4+\ve^{-1}/\log 2$ (and thus $t<T+1<3\ve^{-1}$) since, by the proof of Theorem 2, if $\lfloor t\rfloor\ge 4+\ve^{-1}/\log 2\ge 14$, for $m=\lfloor t\rfloor-4\ge \text{max}\{9,\ve^{-1}/\log 2\}$, we have $D(t-1)\ge HD(C_m)>1-\frac1{m\log 2}>1-\ve$ (and so $D(t-1)>D(t)-\ve$).

Under these hypothesis, we will apply the conclusions of Lemma 2 for $\eta=\ve$. In its proof, in this case, it is enough to assume $r_0\ge 1/\tau^2$ and that, for $k=8N(t,r_0)^2\lceil 2/\tau\rceil$, $\frac{\log N(t,r_0)}{r_0}<(1+\tau/2)\frac{\log N(t,k(r_0-1))}{k(r_0-1)}$ (indeed, assuming the above bounds for $t$, it is not difficult to check that, except for this inequality relating $N(t,r_0)$ and $N(t,k(r_0-1))$, the claims in other parts of the proof of the Lemma that use the assumptions that $r_0$ and $k$ are large are satisfied provided $r_0\ge 1/\tau^2$). 

We define a sequence $(c_n)_{n\ge 0}$ recursively by $c_0=\lceil \frac1{\tau^2}\rceil$ and, for every $n\ge 0$, $c_{n+1}=8N(t,c_n)^2\lceil \frac2{\tau}\rceil(c_n-1)$. We claim that, for some integer $s_0<(1+\frac2{\tau})\log(4/\ve)$, we will have $\frac{\log N(t,c_{s_0})}{c_{s_0}}<(1+\frac{\tau}2)\frac{\log N(t,c_{s_0+1})}{c_{s_0+1}}=(1+\frac{\tau}2)\frac{\log N(t,k(c_{s_0}-1))}{k(c_{s_0}-1)}$, with $k=8N(t,c_{s_0})^2\lceil \frac2{\tau}\rceil$. Indeed, if it is not the case, then $\frac{\log N(t,c_{n+1})}{c_{n+1}}\le(1+\frac{\tau}2)^{-1}\frac{\log N(t,c_n)}{c_n}$ for $0\le n<(1+\frac2{\tau})\log(4/\ve)$, and so, for $M=\lceil(1+\frac2{\tau})\log(4/\ve)\rceil$, we would have 
$\frac{\log N(t,c_M)}{c_M}\le (1+\frac{\tau}2)^{-M}\cdot\frac{\log N(t,c_0)}{c_0}<(\ve/4)\cdot\frac{\log N(t,c_0)}{c_0}$, since $(1+\frac{\tau}2)^{-(1+\frac2{\tau})}<e^{-1}$. On the other hand, it follows by Lemma A3 that, for every $m\ge c_0$, $N(t,m)<(T+1)^2e^m<e^{2m}$ (recall that $c_0=\lceil \frac1{\tau^2}\rceil=\lceil \frac{1600}{\ve^2}\rceil$), and so $\frac{\log N(t,c_M)}{c_M}\le(\ve/4)\cdot\frac{\log N(t,c_0)}{c_0}<\ve/2$. This leads to a contradiction since, for every positive integer $m$, $\ve\le D(t)\le \frac{\log(T^2N(t,m))}m$ and, in particular, $\frac{\log N(t,c_M)}{c_M}>\ve-\frac{2\log T}{c_M}\ge\ve/2$, since $c_M\ge c_0\ge \frac{1600}{\ve^2}$ and $T<3/\ve$. 

Now let $r_0=c_{s_0}$. By the previous discussion, the proof of Lemma 2 works for this $r_0$ (and $k=8N(t,r_0)^2\lceil 2/\tau\rceil$), so, for 
$$\delta=\frac1{(T+1)^{2(m_1+m_2+\tilde m)}}\ge \frac1{(T+1)^{2k\cdot\text{max}\{|\beta|, \beta\in C(t,r_0)\}}}\ge\frac1{(T+1)^{2k\cdot r_0/\log 2}},$$
we have $D(t-\delta)>(1-\ve)D(t)>D(t)-\ve$. We will now give an explicit positive lower bound for $\delta$ in terms of $\ve$. In order to do this we define recursively, for each integer $n\ge 0$ and $x\in\re$, the functions ${\cal T}(x,n)$ and ${\cal T}(n)$ by ${\cal T}(x,0)=x$, ${\cal T}(x,n+1)=e^{{\cal T}(x,n)}$ and ${\cal T}(n)={\cal T}(1,n)$. We have, for every $n\ge 0$, 
$$c_{n+1}=8N(t,c_n)^2\lceil \frac2{\tau}\rceil(c_n-1)\le 8e^{4c_n}\cdot\frac3{\tau}\cdot c_n<e^{e^{c_n}},$$
since, for every $n\ge 0$, $c_n\ge c_0\ge\frac1{\tau^2}=\frac{1600}{\ve^2}$ and $N(t,c_n)\le e^{2c_n}$, therefore $r_0=c_{s_0}<{\cal T}(c_0,2s_0)={\cal T}(\lceil \frac1{\tau^2}\rceil,2s_0)$ and 
$$2\log(T+1)\cdot k\cdot \frac{r_0}{\log 2}=16\log(T+1)\cdot N(t,r_0)^2\lceil \frac2{\tau}\rceil \cdot\frac{r_0}{\log 2}\le 16\log(3/\ve)\cdot e^{4r_0}\cdot\frac3{\tau}\cdot\frac{r_0}{\log 2}<e^{e^{r_0}},$$
so 
$$\delta\ge\frac1{(T+1)^{2k\cdot r_0/\log 2}}=e^{-2\log(T+1)\cdot k\cdot r_0/\log 2}>e^{-e^{e^{r_0}}}>\frac1{{\cal T}(c_0,2s_0+3)}.$$
Finally, since $2^k\ge k^2$ for every $k\ge 4$, it follows by induction that, for every $n\ge 4$, ${\cal T}(n)\ge(n+1)^6$ for every $n\ge 0$. Indeed, ${\cal T}(4)>2^{16}>5^6$ and, for $n\ge 4$, ${\cal T}(n+1)>2^{{\cal T}(n)}\ge {\cal T}(n)^2\ge(n+1)^{12}>(n+2)^6$. This implies that ${\cal T}(\lfloor 1/\ve \rfloor)\ge(1/\ve)^6>1601/\ve^2>\lceil \frac{1600}{\ve^2}\rceil=\lceil \frac1{\tau^2}\rceil=c_0$ (recall that $0<\ve<1/7$), so ${\cal T}(c_0,2s_0+3)<{\cal T}({\cal T}(\lfloor 1/\ve \rfloor),2s_0+3)={\cal T}(\lfloor 1/\ve \rfloor+2s_0+3)$, and, since $s_0<(1+\frac2{\tau})\log(4/\ve)$, we have
$$\lfloor 1/\ve \rfloor+2s_0+3<3+\lfloor 1/\ve \rfloor+2(1+\frac2{\tau})\log(4/\ve)\le 3+1/\ve+2(1+\frac2{\tau})\log(4/\ve)=$$
$$=3+1/\ve+2(1+\frac{80}{\ve})\log(4/\ve)<\frac{161}{\ve}\log(4/\ve),$$
and therefore
$$\delta>\frac1{{\cal T}(c_0,2s_0+3)}>\frac1{{\cal T}(\lfloor 1/\ve \rfloor+2s_0+3)}\ge \frac1{{\cal T}(\lfloor \frac{161}{\ve}\log(4/\ve)\rfloor)}.$$

\section{Appendix: Basic facts and estimates on continued fractions}

We will prove here some elementary facts on continued fractions used in the previous sections. We refer to [CF] for the facts used but not proved here.
 
Let $x = [a_0;a_1,a_2,a_3,\dots]$ be a real number, and $(\frac{p_n}{q_n})_{n \in \mathbb N}, \frac{p_n}{q_n}= [a_0;a_1,a_2,\dots,a_n]$ its sequence of convergents.

  We have, for every $n\in \nb$, $p_{n+1}q_n-p_nq_{n+1}=(-1)^n$,
	$$x=[a_0; a_1, a_2, \dots, a_n, \al_{n+1}]=\frac{\alpha_{n+1}p_n+p_{n-1}}{\alpha_{n+1}q_n+q_{n-1}},$$
	and so
	$$\alpha_{n+1}=\frac{p_{n-1}-q_{n-1}x}{q_{n}x-p_{n}}$$
	and
  $$x-\frac{p_n}{q_n}
  = \frac{(-1)^n}{(\alpha_{n+1}+\beta_{n+1})q_n^2},
  $$
  where
  $$\beta_{n+1}
  = \frac{q_{n-1}}{q_n}
  = [0;a_n,a_{n-1},a_{n-2},\dots,a_1].
  $$
In particular,
\[ \left| x-\frac{p_n}{q_n} \right|= 
\frac 1{(\alpha_{n+1}+\beta_{n+1})q_n^2}. \]

Recall that, given a finite sequence $\al=(a_1,a_2,\dots,a_n)\in(\nb^*)^n$, we define its {\it size\/} by $s(\al):=|I(\al)|$, where $I(\al)$ is the interval $\{x\in[0,1] \mid x=[0; a_1,a_2,\dots,a_n,\al_{n+1}], \al_{n+1}\ge 1\}$, whose endpoints are $p_n/q_n$ and $\frac{p_n+p_{n-1}}{q_n+q_{n-1}}$, $r(\al)=\lfloor\log s(\al)^{-1} \rfloor$ and, for $r \in \nb, P_r= \{\al=(a_1,a_2,\dots,a_n) \mid r(\al) \ge r$ and $r((a_1,a_2,\dots,a_{n-1}))<r\}$.

Since $\alpha_{n+1}=\frac{p_{n-1}-q_{n-1}x}{q_{n}x-p_{n}}$, the $n$-th iterate of the Gauss map restricted to the interval $I(\al)$ is given by
$${g^n}|_{{I(\al)}}(x)=g^n(\alpha_1^{-1})=\alpha_{n+1}^{-1}=\frac{q_{n}x-p_{n}}{p_{n-1}-q_{n-1}x}.$$

\vskip .05in
\noindent
{\bf Lemma A1.} {\it If $a_0, a_1, a_2, \ldots$, $b_n, b_{n+1}, \ldots$ are positive integers with $a_{n-1}, a_n, a_{n+1}, b_n, b_{n+1} \le T$ and $a_n\neq b_n$, then
$$|[a_0;a_1, a_2,\ldots,a_n, a_{n+1}, a_{n+2}, \ldots]-[a_0;a_1, a_2,\ldots,b_n, b_{n+1}, b_{n+2}, \ldots]|>\frac1{(T+1)(T+2)q_n^2},$$
where $q_n=q_n(a_1, a_2,\ldots, a_n)$.}

\noindent
{\bf Proof:} Let $x=[a_0;a_1, a_2,\ldots,a_n, a_{n+1}, a_{n+2}, \ldots]$ and $y=[a_0;a_1, a_2,\ldots,a_{n-1}, b_n, b_{n+1}, b_{n+2}, \ldots]$. We have 
$$x=\frac{\alpha_{n}(x)p_{n-1}+p_{n-2}}{\alpha_{n}(x)q_{n-1}+q_{n-2}}, y=\frac{\alpha_{n}(y)p_{n-1}+p_{n-2}}{\alpha_{n}(y)q_{n-1}+q_{n-2}},$$
and so
$$|x-y|=\left|\frac{\alpha_{n}(x)p_{n-1}+p_{n-2}}{\alpha_{n}(x)q_{n-1}+q_{n-2}}-\frac{\alpha_{n}(y)p_{n-1}+p_{n-2}}{\alpha_{n}(y)q_{n-1}+q_{n-2}}\right|=\left|\frac{(\alpha_{n}(x)-\alpha_{n}(y))(p_{n-1}q_{n-2}-p_{n-2}q_{n-1})}{(\alpha_{n}(x)q_{n-1}+q_{n-2})(\alpha_{n}(y)q_{n-1}+q_{n-2})}\right|$$
$$=\left|\frac{(\alpha_{n}(x)-\alpha_{n}(y))(-1)^{n-2}}{(\alpha_{n}(x)q_{n-1}+q_{n-2})(\alpha_{n}(y)q_{n-1}+q_{n-2})}\right|=\frac{|\alpha_{n}(x)-\alpha_{n}(y)|}{(\alpha_{n}(x)q_{n-1}+q_{n-2})(\alpha_{n}(y)q_{n-1}+q_{n-2})}.$$
We have $\lfloor \alpha_n(x) \rfloor=a_n$, $\lfloor \alpha_n(y) \rfloor=b_n$ and $a_n\neq b_n$, so \hfill\break $|\alpha_{n}(x)-\alpha_{n}(y)|>[1;T+1]-[0;1,T+1]=\frac1{T+1}+\frac1{T+2}>\frac2{T+2}$. Moreover, $\alpha_{n}(x)q_{n-1}+q_{n-2}<(1+a_n)q_{n-1}+q_{n-2}=q_n+q_{n-1}<2q_n$ and $\alpha_{n}(y)q_{n-1}+q_{n-2}<\alpha_{n}(y)(q_{n-1}+q_{n-2})\le \alpha_{n}(y)q_n<(T+1)q_n$, and so 
$$|x-y|>\frac2{T+2} \cdot \frac1{2q_n\cdot (T+1)q_n}=\frac1{(T+1)(T+2)q_n^2}.\qed$$

\vskip .05in
\noindent
{\bf Lemma A2.} {\it If $\alpha=a_1a_2\cdots a_m$ and $\beta=b_1b_2\cdots b_n$ are finite words, then
$$\dfrac{1}{2}s(\alpha)s(\beta)<s(\alpha\beta)<2s(\alpha)s(\beta).$$}
\noindent
{\bf Proof:} By Euler's property of continuants (cf. Appendix 2 of [CF]), we have
$$q_{m+n}(\alpha\beta)=q_m(\alpha)q_n(\beta)+q_{m-1}(a_1a_2\cdots a_{m-1})q_{n-1}(b_2b_3\cdots b_n)$$
and then
$$q_m(\alpha)q_n(\beta)<q_{m+n}(\alpha\beta)<2q_m(\alpha)q_n(\beta).$$
From the left inequality, we have
\begin{eqnarray*}
s(\alpha\beta)&=&\dfrac{1}{q_{m+n}(\alpha\beta)\left[q_{m+n}(\alpha\beta)+q_{m+n-1}(\alpha\beta)\right]}\\
              &<&\dfrac{1}{q_m(\alpha)q_n(\beta)\left[q_m(\alpha)q_n(\beta)+q_m(\alpha)q_{n-1}(\beta)\right]}\\
&<&\dfrac{2}{q_m(\alpha)q_n(\beta)\left[q_m(\alpha)+q_{m-1}(\alpha)\right]\left[q_n(\beta)+q_{n-1}(\beta)\right]}\\ &=&2\cdot\dfrac{1}{q_m(\alpha)\left[q_m(\alpha)+q_{m-1}(\alpha)\right]}\cdot\dfrac{1}{q_n(\beta)\left[q_n(\beta)+q_{n-1}(\beta)\right]}\\
&=&2s(\alpha)s(\beta),
\end{eqnarray*}
where in the second inequality we used that
\begin{eqnarray*}
2q_m(\alpha)q_n(\beta)+2q_m(\alpha)q_{n-1}(\beta)&>&q_m(\alpha)q_n(\beta)+q_m(\alpha)q_{n-1}(\beta)+\\
&&q_{m-1}(\alpha)q_n(\beta)+q_{m-1}(\alpha)q_{n-1}(\beta)\\
\iff\hspace{.5cm} q_m(\alpha)q_n(\beta)+q_m(\alpha)q_{n-1}(\beta)&>&q_{m-1}(\alpha)q_n(\beta)+q_{m-1}(\alpha)q_{n-1}(\beta),
\end{eqnarray*}
which is obviously true. For the other inequality, proceed analogously:
\begin{eqnarray*}
s(\alpha\beta)&=&\dfrac{1}{q_{m+n}(\alpha\beta)\left[q_{m+n}(\alpha\beta)+q_{m+n-1}(\alpha\beta)\right]}\\
&>&\dfrac{1}{2}\cdot\dfrac{1}{q_m(\alpha)q_n(\beta)\left[q_{m+n}(\alpha\beta)+q_{m+n-1}(\alpha\beta)\right]}\\
&>&\dfrac{1}{2}\cdot\dfrac{1}{q_m(\alpha)q_n(\beta)\left[q_m(\alpha)+q_{m-1}(\alpha)\right]\left[q_n(\beta)+q_{n-1}(\beta)\right]}\,,
\end{eqnarray*}
as 
\begin{eqnarray*}
q_{m+n}(\alpha\beta)+q_{m+n-1}(\alpha\beta)&<&\left[q_m(\alpha)+q_{m-1}(\alpha)\right]\left[q_n(\beta)+q_{n-1}(\beta)\right]\\
\iff\hspace{.2cm} q_{m-1}(\alpha){\tilde q}_{n-1}(\beta)+q_{m-1}(\alpha){\tilde q}_{n-2}(\beta)&<&
q_{m-1}(\alpha)q_n(\beta)+q_{m-1}(\alpha)q_{n-1}(\beta)\\
\iff\hspace{2.75cm}{\tilde q}_{n-1}(\beta)+{\tilde q}_{n-2}(\beta)&<&q_n(\beta)+q_{n-1}(\beta),
\end{eqnarray*}
where ${\tilde q}_{n-1}(\beta)=q_{n-1}(b_2b_3\cdots b_n)$ and ${\tilde q}_{n-2}(\beta)=q_{n-2}(b_2b_3\cdots b_{n-1})$, and the last inequality is true, since we have ${\tilde q}_{n-1}(\beta)<q_n(\beta)$ and ${\tilde q}_{n-2}(\beta)<q_{n-1}(\beta)$. This concludes the proof. \qed

\noindent
{\bf Lemma A3.} If $\al=(a_1,a_2,\dots,a_n)\in(\nb^*)^n$ belongs to $P_r$ and $1\le a_j\le T$ for $1\le j\le n$, then $s(\al)>((T+1)^2e^r)^{-1}$.

\noindent
{\bf Proof:} We have $r(a_1,a_2,\dots,a_{n-1})<r$, so $s(a_1,a_2,\dots,a_{n-1})=(q_{n-1}(q_{n-1}+q_{n-2}))^{-1}>e^{-r}$, and thus $s(\al)^{-1}=q_n(q_n+q_{n-1})\le (Tq_{n-1}+q_{n-2})((T+1)q_{n-1}+q_{n-2})<(T+1)q_{n-1}\cdot (T+1)(q_{n-1}+q_{n-2})=(T+1)^2q_{n-1}(q_{n-1}+q_{n-2})<(T+1)^2e^r$, so $s(\al)>((T+1)^2e^r)^{-1}$. \qed
\vskip .1in
\noindent
{\bf Proposition A1.} 
{\it Let $r, k$ be positive integers and $\alpha_i, 1\le i \le k$ finite sequences which belong to $P_r$ and whose elements are bounded by $T$. Then, if $\al=\al_1\al_2\cdots\al_k$, we have $s(\alpha)>(2(T+1)^2 e^r)^{-k}$.}

\noindent
{\bf Proof:} For $1\le i\le k$ we have, by lemma A3, $s(\al_i)>((T+1)^2e^r)^{-1}$. So, using the lemma A2, we get
$$s(\al)=s(\al_1\al_2\cdots\al_k)>\left(\frac12\right)^{k-1}s(\al_1)s(\al_2)\cdots s(\al_k)>$$
$$>\left(\frac12\right)^{k-1}(((T+1)^2e^r)^{-1})^k>(2(T+1)^2 e^r)^{-k}.\,\,\,\,\,\,\,\,\,\,\,\,\,\,\,\, \qed$$ 
\vskip .2in
\noindent{\bf References}

\vglue .1in

[B] Bugeaud, Y. - Sets of exact approximation order by rational numbers. II - Unif. Distrib. Theory v. 3 , n. 2, p. 9-20, 2008.

[CF] Cusick, T. and Flahive, M. - The Markov and Lagrange spectra -
Mathematical Surveys and Monographs, n. 30 - AMS.

[Fa] Falconer, K. -The geometry of fractal sets - Cambridge Tracts in Mathematics, 85. Cambridge University Press, Cambridge, 1986.

[FMM] Ferenczi, S., Mauduit, C. and Moreira, C.G. - An algorithm for the word entropy - https://arxiv.org/abs/1803.05533

[F] Freiman, G.A. - Diofantovy priblizheniya i geometriya chisel (zadacha Markova) [Diophantine approximation and geometry of numbers (the Markov spectrum)], Kalininskii Gosudarstvennyi Universitet, Kalinin, 1975.

[H] Hall, M. - On the sum and products of continued fractions - Annals of Math. v. 48, p. 966-993, 1947.

[J] Jarn\'\i k, V. - Zur metrischen Theorie der diophantischen Approximationen - Prace Mat.-Fiz. 36, p. 91-106, 1928.

[Ma] Markov, A. - Sur les formes quadratiques binaires indéfinies, Math. Ann. v. 15, p. 381-406, 1879.

[MM1] Matheus, C. and Moreira, C.G., \ $HD(M\setminus L) > 0.353$ -  https://arxiv.org/abs/1703.04302

[MM2] Matheus, C. and Moreira, C.G., \ $HD(M\setminus L) < 0.986927$ - https://arxiv.org/abs/1708.06258

[MM3] Matheus, C. and Moreira, C.G., \ Markov spectrum near Freiman's isolated points in $M\setminus L$ - https://arxiv.org/abs/1802.02454

[MM4] Matheus, C. and Moreira, C.G., \ New numbers in $M\setminus L$ beyond $\sqrt{12}$: solution to a conjecture of Cusick. https://arxiv.org/abs/1803.01230

[Mo] Moreira, C.G. - Geometric properties of images of cartesian products of regular Cantor sets by differentiable real maps. https://arxiv.org/abs/1611.00933

[MY] Moreira, C.G. and Yoccoz, J.-C. - Stable intersections of regular Cantor
sets with large Hausdorff Dimensions - Annals of Math. v. 154, n.1, p. 45 - 96, 2001.

[P] Perron, O. - \"Uber die Approximation irrationaler Zahlen durch rationale II - S.-B. Heidelberg Akad. Wiss., Abh. 8, 1921, 12 pp.

[PT] Palis, J. and Takens, F. - Hyperbolicity and sensitive chaotic dynamics at homoclinic bifurcations, Fractal dimensions and infinitely many attractors. Cambridge Studies in Advanced Mathematics, 35. Cambridge University Press, Cambridge, 1993. x+234 pp.

\end{document}